%

\input amstex
\documentstyle{amsppt}

   \NoBlackBoxes 
 
\loadmsbm

\def\p#1{\, #1}
\def\rmp#1{\,{\rm #1}}

\def\comma{{\rm ,}\ }
\def\({\/\ {\rm (}}
\def\){\/{\rm )\ }}
\def\semicolon{\/{\rm ;}\ }
\def\colon{\/{\rm :}\ }

\def\fSg#1#2{\frak {Sg}^{\frak{#1}}#2}
\def\fSgx#1{\frak {Sg}@,@,@,#1}

\def\Sg#1#2{S@!@!@!@!@!g^{\frak{#1}}#2}
\def\Sgx#1{S@!@!@!@!@!g@,@,@,#1}

\def\fA{{\frak A}}
\def\fB{{\frak B}}
\def\fC{{\frak C}}

\def\fF{{\frak F}}
\def\fM{{\frak M}}
\def\fN{{\frak N}}
\def\fP{{\frak P}}
\def\fU{{\frak U}}
\def\fV{{\frak V}}
\def\fW{{\frak W}}
\def\fZ{{\frak Z}}

\def\ka{\kappa}
\def\EM#1{EM@,({\frak #1},\Phi)}
\def\EMM#1{EM@,({\frak #1},\Psi)}
\def\prc{\preccurlyeq}
\def\nprc{\not\preccurlyeq}
\def\rst#1#2{{#1}\!\!\restriction\!@!@!@!{#2}}
\def\pe{p\!\restriction\!{\!}^*\!E}
\def\lr#1#2{\langle\,{#1}\,,\,{#2}\,\rangle}
\def\fpr#1{\frak {Pr}_a{({#1})}}
\def\diff{-}
\def\preom#1{{ }^{\omega\sssize\geq}#1}
\def\num#1#2{\medskip\smallskip\flushpar (#1)\quad #2}
\def\midnum#1#2{\medskip\flushpar (#1)\quad #2}
\def\endnum{\medskip\smallskip\flushpar}
 
\def\*{\!*@!@!@!@!@!}
\def\jx{J\*X_\alpha}
\def\ky{K\*Y_\alpha}
\def\mybigcup{{\tsize\bigcup}}

\document
\topmatter
\title Universal Theories Categorical in Power and $\kappa$-generated Models
\endtitle \author Steven Givant and Saharon Shelah\endauthor
\thanks Shelah's research was partially supported by the National Science
Foundation and the United States-Israel Binational Science Foundation.  This
article is item number 404 in Shelah's bibliography. The authors
would like to thank Garvin Melles for reading a preliminary draft of
the paper and making several very helpful suggestions. \endthanks \address
Mills College\endaddress \address Hebrew University,
Rutgers University, and
MSRI
\endaddress
\abstract  We
investigate a notion called {\it uniqueness in power $\ka$} that is
akin to categoricity in power $\ka$, but is based on the  cardinality
of the generating sets of models instead of on the cardinality of their
universes.  The notion is quite useful for formulating
categoricity-like questions regarding powers below the cardinality of
a theory. We prove, for (uncountable) universal theories $T$, that if
$T$ is $\kappa$-unique for one uncountable $\kappa$, then it is $\kappa$-unique for every
uncountable $\kappa$; in particular, it is categorical in powers greater than
the cardinality of $T$. \endabstract  \endtopmatter

It is well known that the notion of categoricity in power exhibits certain
irregularities in ``small" cardinals, even when applied to such simple theories as
universal Horn theories.  For example, a  countable universal Horn theory categorical in
one uncountable power is necessarily categorical in all  uncountable powers, by
Morley's theorem, but it need not be countably categorical.

Tarski suggested that, for universal Horn theories $T$, this irregularity
might be overcome by replacing the notion of categoricity in power by that of
freeness in power. $T$ is {\it free in power\/} $\ka$, or $\ka$-{\it free\/}, if
it has a model of power $\ka$ and if all such models are free, in the general
algebraic sense of the word, over the class of all models of $T$.  $T$ is a
{\it free theory\/} if each of its models is free over the class of all models
of $T$.  It is trivial to check that, for $\ka>|T|$, categoricity and
freeness in power $\ka$ are the same thing. For $\ka\leq |T|$ they are
not the same thing.  For example, the (equationally axiomatizable) theory of
vector spaces over the rationals is an example of a free theory, 
categorical in every uncountable power, that is $\omega$-free, but not 
$\omega$-categorical. Tarski formulated the following problem: Is a  universal
Horn theory that is free in one infinite power necessarily free  in all
infinite powers? Is it a free theory?

Baldwin, Lachlan, and McKenzie, in Baldwin-Lachlan \cite{1973}, and
Pal\-yu\-tin,  in Aba\-ku\-mov-Pal\-yu\-tin-Shish\-ma\-rev-Taits\-lin
\cite{1973}, proved that a countable $\omega$-categorical universal Horn theory
is $\omega_1$-categorical, and hence categorical in all infinite powers.
Thus, it is free in all infinite powers.  Givant \cite{1979} showed that
a countable $\omega$-free, but not $\omega$-categorical, universal Horn theory
is  also $\omega_1$-categorical, and in fact it is a free theory.   Further, he
proved that a universal Horn theory, of any cardinality  $\ka$, that is
$\ka$-free, but not $\ka$-categorical, is necessarily a  free theory.

Independently, Baldwin-Lachlan, Givant, and Palyutin all found
examples of countable $\omega_1$-categorical universal Horn theories that
are not $\omega$-free, and of countable $\omega_1$- and
$\omega$-categorical universal Horn theories that are not free theories. 

Thus, Tarski's implicit problem remains: find a notion akin to
categoricity in power that is regular for universal Horn theories, i.e.,
if it holds in one infinite power, then it holds in every infinite power.

One of the difficulties with the notions of categoricity in power and
freeness in power is that they are defined in terms of the cardinality of
the universes of models instead of the cardinality of the generating sets.  This
causes difficulties when trying to work with  powers $<|T|$.  

Let's call a model $\fA$ {\it
strictly\/} $\ka$-{\it generated\/} if $\ka$ is the minimum of the cardinalities of 
generating sets of $\fA$.  We define a theory $T$ to be $\ka$-{\it
unique\/} if it has, up to isomorphisms, exactly one strictly $\ka$-generated model. 
For cardinals $\ka > |T|$, the  notions of $\ka$-categoricity, $\ka$-freeness, and
$\ka$-uniqueness coincide (in the case when $T$ is universal Horn). When $\ka=|T|$, we
have $$
\ka@-categoricity\quad\Rightarrow\quad\ka@-freeness\quad\Rightarrow
\quad\ka@-uniquenes\quad,$$
but none of the reverse implications hold.

Givant \cite{1979}, p\. 24, asked, for universal Horn theories $T$, whether
$\ka$-uniqueness is the regular notion that Tarski was looking for, i.e.,
(1) Does $\ka$-uniqueness for one infinite $\ka$ imply it for all infinite $\ka$?
For countable $T$, he answered this question affirmatively by showing that
$\omega$-uniqueness is equivalent to categoricity in uncountable powers.  For
uncountable $T$, he provided a partial affirmative answer by showing
that categoricity in power $>|T|$ implies $\ka$-uniqueness for all
infinite $\ka$, and is, in turn, implied by $\ka$-uniqueness when $\ka=|T|$ and
$\ka$ is regular.  However, the problem whether $\kappa$-uniqueness
implies categoricity in powers $> |T|$ when $\omega\leq\ka<|T|$, or when $\omega
<\ka = |T|$ and $\ka$  is singular, was left open.

In this paper we shall prove the following:

\proclaim{Theorem} 
A universal theory $T$ that is
$\ka$-unique for some $\ka>\omega$ is $\ka$-unique for every $\ka >\omega$. 
In particular\comma it is categorical in powers $>|T|$.
\endproclaim

It follows from the previous remarks that a universal Horn theory $T$ which is
$\ka$-unique for some $\ka>\omega$ is $\ka$-unique for every $\ka\ge\omega$. 
Thus, the only part of (1) that still remains open is the case 
when $T$ is uncountable and $\omega$-unique. A more general formulation of 
this open problem is the following:

\proclaim{Problem}  
Is an $\omega$-unique universal theory $T$ necessarily
categorical in powers $>|T|$\rmp? In particular\comma is a countable
$\omega$-unique \(or $\omega$-categorical\,\) universal theory
$\omega_1$-categorical\rmp?\endproclaim

An example due to Palyutin,
in  Abakumov-Palyutin-Shishmarev-Taitslin \cite{1973}, shows that a countable
universal theory categorical in uncountable powers need not be
$\omega$-unique.  In fact, in
Palyutin's example the finitely
generated models are all finite, and there are countably many non-isomorphic,
strictly  $\omega$-generated models.  Thus, for universal theories, $\ka$-uniqueness for
some $\ka>\omega$ does not imply $\omega$-uniqueness. 

To prove our theorem, we shall show that, under the given 
hypotheses, the theory of the infinite models of $T$ is complete,
superstable, and unidimensional, and that all sufficiently large
 models are $a$-saturated. Thus, we shall make use of some of
the notions and results of stability theory that are developed in Shelah \cite{1990}
(see also Shelah \cite{1978}). We will assume that the reader is acquainted with the
elements of model theory and with such basic notions from stability theory as
superstability, $a$-saturatedness, strong type, regular type, and Morley
sequence.  We begin by reviewing some notation and terminology, and then
proving a few elementary lemmas.

The letters $m$ and $n$ shall denote finite cardinals, and $\kappa$
and $\lambda$ infinite cardinals. The cardinality of a set $U$ is
denoted by $|U|$.  The set-theoretic
difference of $A$ and $B$ is denoted by $A\diff B$.  If $\vartheta$ is a function, and
$\bar a = \langle a_0,\dots, a_{n-1}\rangle$ is a sequence of elements in the
domain of $\vartheta$, then $\vartheta (\bar a)$ denotes the sequence $\langle
\vartheta(a_0),\dots, 
 \vartheta(a_{n-1})\rangle$.
We denote the restriction of $\vartheta$ to a subset $X$ of its domain by
$\rst \vartheta X$, and a similar notation is employed for the
restriction of a relation.  A sequence $\langle X_\xi: \xi < \lambda
\rangle$ of sets is {\it increasing\/} if $X_\xi\subseteq X_\eta$
for $\xi < \eta <\lambda$, and {\it continuous\/} if $X_\delta = \bigcup_{\xi <
\delta} X_\xi$ for limit ordinals $\delta<\lambda$.

 We use German letters $\fA, \fB, \fC,...$ to denote
models, and the  corresponding Roman letters $A,B,C,...$ to denote their
respective  universes. If $\tau(x_0,...,x_{n-1})$ is a term in a (fixed) 
language $L$ for  $\fA$, and if and $\bar a$ an $n$-termed sequence of elements
in $A$,  symbolically $\bar a\in{}^n A$, then the value of $\tau$  at $\bar
a \text{ in } \fA$ is denoted by $\tau^{\fA}[\bar a],$ or simply by  $\tau[\bar
a]$. A similar notation is used for formulas.  Suppose $\bar a\in{}^n A$
and $X\subseteq A$.  The type of  $\bar a$ over $X$ (in $\fA$), i.e. the set of
formulas in the language of  ${\lr\fA x}_{x \in X}$ that are satisfied by $\bar
a$ in  the latter model, is denoted by $tp^{\fA}(\bar a,X)$, or simply by $tp(\bar
a,X)$, when no confusion can arise.  The strong type of
$\bar a$ over $X$ (in $\fA$), i.e., the set of  formulas in the language of ${\lr\fA
b}_{b\in A}$ that are almost over $X$ and that are satisfied by  $\bar a$, is
denoted by $stp^{\fA}(\bar a,X)$, or simply by $stp(\bar a,X)$. 
If $p(\bar x)$ is a strong type and
$E \subseteq A$ is a base for $p$ in $\fA,$ then $\pe$ denotes the set of
formulas in $p$ that are almost over $E$.

We write $\fA \subseteq \fB$ to express that $\fA$ is a submodel of $\fB$. The
submodel of $\fA$ generated by a set $X \subseteq A$ is denoted by 
$\fSg A {(X)}$, or simply by $\fSgx {(X)}$, and its universe by $\Sgx 
{(X)}$. A model is $\mu$-{\it generated\/} if it is generated by a set of
cardinality $\mu$, and {\it strictly} $\mu$-{\it generated\/} if it is
$\mu$-generated, but not $\nu$-generated for any $\nu<\mu$. Every model $\fA$ has a
generating set of minimal cardinality, and hence is strictly $\mu$-generated for some
(finite or infinite) $\mu$.  If $X$ is  a generating set of $\fA$ of minimal
cardinality, and $Y$ is any other generating set of $\fA$, then  there must be
a subset $Z$ of $Y$ of power at most $|X| + \omega$ such that $Z$ generates $X$,
and hence also $\fA$.
A set
$X \subseteq A$ is {\it irredundant\/} (in $\fA$) if, for every
$Y\subsetneqq X$ we  have $\fSgx {(Y)} \neq \fSgx {(X)}.$ 

 A model $\fA$ is
an 
 $n$-{\it submodel\/} of $\fB$, and $\fB$ an
$n$-{\it extension\/} of $\fA$, in symbols $\fA \prc_n \fB$, if for every
$\Sigma_n\,$-formula $\varphi(x_0,...,x_{k-1})$ in the language of $\fA$, and
every  $\bar a\in{}^k A$, we have $\fA \models\varphi[\bar a]$
iff $\fB \models \varphi[\bar a]$.  A $0$-submodel of $\fB$ is just a submodel in
the usual sense of the word, and an elementary submodel---in symbols,
$\fA\prc\fB$---is just an $n$-submodel for each $n$.  We write $\fA\prec\fB$ to
express the fact that $\fA$ is a {\it proper\/} elementary submodel of
$\fB$, i.e., $\fA\prc\fB$ and $\fA\neq\fB$. A theory is {\it model complete\/}
if, for any two models $\fA$ and $\fB$, we have $\fA\prc\fB$ iff
$\fA\subseteq\fB$.  For any theory $T$, we denote by $T_\infty$ the theory of
the infinite models of $T$.

We shall always use the phrase {\it dense ordering\/} to mean a non-empty
dense linear ordering without endpoints.  It is well-known that the theory of
such orderings admits elimination of quantifiers.  Hence, in any such ordering
$\fU =\,\lr U <$, if $\bar a$ and $\bar b$ are two sequences in ${}^n U$
that are {\it atomically\/} {\it equivalent\/}, i.e., that satisfy the same
atomic formulas, then they are {\it elementarily equivalent\/}, i.e., they
satisfy the same elementary formulas.

Fix a linear ordering $\fV = \lr V <$, and set
$W=\preom V = \bigcup_{\mu\leq\omega} {}^\mu V $. Let $\lessdot$ be the
(proper) initial segment relation on $W$, and, for each $\mu\leq\omega$, let
$P_\mu$ be the set of elements in $W$ with domain $\mu$, i.e., $P_\mu = {}^\mu
V$.  There is a natural lexicographic ordering, $<$\ , on $W$ induced by the
ordering of $V$: $f<g$ iff either $f\lessdot g$ or else there is a natural
number $n$ in the domain of both $f$ and $g$ such that $\rst f n = \rst g n$ and
$f(n) < g(n)$ in $\fV$. Take $h$ to be the binary function on $W$ such
that, for any $f,g$ in $W$, $h(f,g)$ is the greatest common initial segment
of $f$ and $g$.  We shall call the structure
$\fW=\langle\,W\,,\,<\,,\,\lessdot\,,\,P_\mu\,,\,h\,\rangle_{\mu
\leq\omega}$
the {\it full tree structure  over
$\fV$
with $\omega + 1$ levels{\rm, or, for short, the} full tree  over
$\fV$.} 
Any substructure of the full tree over  $\fV$  that is {\it downward
closed\/}, i.e., closed under initial segments, is called a 
 {\it tree over\/} $\fV$.  A {\it tree\/} is any structure isomorphic
to a tree structure over some ordering. A tree $\fU$ over a dense order $\fV$ is
itself called {\it dense\/} if: (i) for every $f$ in $U$ with finite domain, say
$n$, the set of immediate successors of $f$ in $U$, i.e., the set of extensions
of $f$ in $U$ with domain $n+1$, is densely ordered by $<$ in $\fU$, or, put a
different way,  $\{g(n): f\lessdot g \}$ is dense under the ordering
inherited from $\fV$; (ii) every element in $U$ with a finite domain is an
initial segment of an element in $U$ with domain $\omega$.  Just as with
dense orderings, the theory of the class of dense trees admits elimination of
quantifiers.

A model $\fA\,$ is $\ka$-{\it homogeneous\/} if, for every
cardinal $\mu < \ka$  and every pair $\bar a,\bar b \in {}^\mu A$ of
elementarily equivalent sequences,  there is an automorphism  $\vartheta$ of
$\fA$ taking $\bar a$ to  $\bar b$. 
It is well-known that, for regular cardinals $\ka$, any model has
$\ka$-homogenous  elementary extensions (usually of large cardinality). It
follows from our remarks  above that, if $\fU$ is a $\ka$-homogeneous dense
ordering or tree, and if  $\bar a, \bar b\in {}^\mu U$ are atomically equivalent
(where $\mu < \ka$), then there is an automorphism of $\fU$ taking $\bar a$ to
$\bar b$.  A model $\fA$ is
{\it weakly $\omega$-homogenous\/} provided that, for every $n$,  every pair of
sequences $\bar a,\bar a'\in {}^n A$ that are elementarily equivalent, and
every $b\in A$, there is a $b'\in A$ such that $\bar a\hat{\ }\langle b\rangle$
and $\bar a'\hat{\ }\langle b'\rangle$ are elementarily equivalent.  It is well
known that a countable, weakly $\omega$-homogenous model is
$\omega$-homogeneous.  Dense orderings and dense trees are always weakly
$\omega$-homogeneous.  

We turn, now, to some notions from stability theory.  A model $\fA$ is $a$-{\it
saturated\/} if, for any strong type $p$ (consistent with the theory of ${\lr\fA
a}_{a\in A}$),  if $E \subseteq A$ is a finite base for $p$, then $\pe$ is realized in
$\fA$. We say that $\fA$ is 
$a$-{\it saturated in\/} $\fB$, in symbols
$\fA\prc_a\fB$, provided that $\fA\prc\fB$ and that, for any  strong type $p$ of
$\fB$ based on a finite subset $E$ of $A$, if  $\pe$ is realized in $\fB$, then
it is realized in $\fA$.  A model $\fA$  is 
$a$-prime over a set $X$ if
(i) $X\subseteq A$ and $\fA$ is $a$-saturated;
(ii) whenever $\fB$ is an $a$-saturated model elementarily equivalent to
$\fA$ and such that $X\subseteq B$, then  there is an elementary embedding of $\fA$
into $\fB$ that leaves the elements of $X$ fixed. (Here, we assume that $\fA$ and $\fB$ 
are elementary substructures of some monster model).  As is shown in Shelah
\cite{1990}, Chapter IV, Theorem 4.18, for complete, superstable theories, the
$a$-prime model over $X$ exists and is unique, up to isomorphic copies over $X$. 
We shall denote it by $\fpr X$. Let $\fA$ be a model and $B, C$ subsets of $A$
with $C\subseteq B$. A type $p(\bar x)$ over $B$ {\it splits over\/}  $C$ if
there are $\bar b, \bar c$ from $B$ such that $tp(\bar b,C) =tp(\bar c,C)$,
and there is a formula $\varphi(\bar x,\bar y)$ over $C$ such that
$\varphi(\bar x,\bar b)$ and $\neg\varphi(\bar x,\bar c)$ are both in $p$.

We now fix a universal theory $T$ in a language $L$ of arbitrary cardinality. 
In what follows let $L'$ be an expansion of $L$ with built-in Skolem
functions, and $T'$  any Skolem theory in $L'$ that extends $T$.  Without
loss of generality, we may assume that, for every term $\tau(x_0,...,x_{n-1})$
of $L'$ there is a function symbol $f$ of $L'$ of rank $n$ such that the
equation $f(x_0,...,x_{n-1})=\tau$ holds in $T'.$  

For every infinite ordering $\fU=\lr U\leq$ and 
(complete) Ehrenfeucht-Mostowski set $\Phi$ of formulas of $L'$
compatible with $T'$,  there is a model
$\fM'=\EM U$ of $T'$---called a ({\it standard\/}) 
{\it Ehrenfeucht-Mostowski model of $T'$\/}---such that $\fM'$ is generated by $U$ (in
particular, $U \subseteq M'$), and $U$ is a set of $\Phi$-indiscernibles  in $\fM'$ with
respect to the ordering of $\fU$, i.e., if $\varphi(x_0,...,x_{n-1})\in\Phi$, and if $\bar
a$ $\in{}^n U$ satisfies $a_i < a_j$ for $i < j< n$, then
$\fM' \models \varphi[\bar a]$. 

Shelah \cite{1990}, Chapter VII, Theorem 3.6, establishes the existence of certain
generalizations of Ehren\-feucht-Mos\-tow\-ski models.  Suppose that $T$ is,
e.g., nonsuperstable, and that the sequence $\langle \varphi_n(\bar
x,\bar y) :n<\omega\rangle$  of formulas witnesses this
nonsuperstability (see, e.g., Shelah \cite{1990}, Chapter II, Theorems 3.9
and 3.14). Then there is a {\it generalized\/} Ehrenfeucht-Mostowski set
$\Phi$ from which we can construct, for every tree $\fU$, a {\it
generalized\/} Ehrenfeucht-Mos\-tow\-ski model $\fM'=\EM U$ of $T'$.
In other words, given a tree $\fU$, there a model $\fM'$ of $T'$
 generated by a sequence $\langle \bar a_u:u\in U\rangle$ of
$\Phi$-indiscernibles with respect to the atomic formulas of $\fU$; in
particular, if $\bar w$ and $\bar v$ in ${}^n U$ are atomically equivalent in
$\fU$, then $\langle \bar a_{w_0},\dots, \bar a_{w_{n-1}}\rangle$ and 
$\langle \bar a_{v_0},\dots, \bar a_{v_{n-1}}\rangle$ are elementarily
equivalent in $\fM'$.  Moreover, for $w$ in $P_n^{\fM'}$ and $v$ in $P_\omega^{\fM'}$,
we have that $\fM'\models \varphi_n[\bar a_w,\bar a_v]$ iff $w\lessdot v$.  In general,
the sequences  $\bar a_u$ may be of length greater than $1$.  However, to
simplify notation we shall act as if they all have length $1$, and in fact,
we shall identify  $\bar a_u$ with $u$.  Thus, we shall assume that $\fM'$ is
generated by $U$ and that two sequences from $U$ which are atomically equivalent
in $\fU$ are elementarily equivalent in $\fM'$.

 Here are some well-known facts about
Ehrenfeucht-Mostowski models.  Let $\fU , \fV$ be infinite orderings or trees,
and $\fM' = \EM U\ , \fN' = \EM V$. (In the case of trees, we  must assume that
$\fM'$ and $\fN'$ exist.)
\redefine\fSg#1#2{\frak {Sg}^{\frak{#1}}\!@!@!@!(#2)}
\redefine\fSgx#1{\frak {Sg}@,@,@,(#1)}

\proclaim{Fact 1}
If $\vartheta$ is any
embedding of $\fV$ into $\fU$\rmp, then the canonical extension
of $\vartheta$ to $N'$ is an elementary embedding of
$\fN'$ into $\fM'$.\endproclaim 

\noindent In particular,
\proclaim{Fact 2}
Any automorphism of $\fU$ extends to an automorphism of $\fM'$.\endproclaim

\proclaim{Fact 3}
If $\fV\subseteq\fU$\rmp, then $\fSg {\fM'} V\prc\fM'$\rmp, and $\fSg {M'} V$
is isomorphic to $\fN'$ via a canonical isomorphism that is the identity on
$V$.\endproclaim 

\proclaim{Fact 4}
If $\,\fU$ is a linear order\comma then $U$ is irredundant in $\fM'$. 
If $\,\fU$ is a tree\comma then for any element $f$  and subset $X$ of $U$\comma
if $f$ is not an initial segment of any element  of $X$\comma then $f$ is not 
generated by $X$ in $\fM'$.
\endproclaim

We shall always denote the reduct of $\fM'$ to $L$ by $\fM$, in symbols,
$\fM=\rst{\fM'}L$. Facts 1, 2, and 4 transfer automatically from $\fM'$
to $\fM$.  However, $\fM$ is usually not generated by $U$.  We shall
therefore formulate versions of the above facts that apply to $\fSg M U$, and,
more generally, to a collection of models that lie between $\fSg M U$ and
$\fSg {M'} U = \fM'$.
 
\proclaim{Definition}
Suppose $\fM'=\EM U$\rmp, and $K\subseteq L'-L$\.
\roster
\item"{(i)}" For each $f\in K$\comma let $R_f$ be the
relation corresponding to $f^{\fM'}$\comma i.e.\comma $R_f[a_0,\ldots,a_n]$ iff 
$f^{\fM}[a_0,\ldots,a_{n-1}]=a_n$.  Set
$\fM_K={\lr\fM {R_f}}_{f\in K}$.
\item"{(ii)}" For each $V\subseteq U$ set
$$K\*V= V\cup\{f^{\fM'}[\bar a]:f\in K , \bar a 
\text{ from } V\}\quad .$$ \endroster
\endproclaim

Thus, the elements of $K\*V$ are just the elements of $V$, together with the
elements of $\fM'$ that can be obtained from  sequences of $V$ by  single
applications of $f$ in $K$.

\proclaim{Lemma 1} Let $\fU\,, \fV$ be dense
orderings or dense trees\comma $\fM'=\EM U\ , \fN'=\EM V$\comma and $K\subseteq L'\diff
L$.  \(Thus\comma in the case of dense trees\comma we postulate that generalized
Ehren\-feucht-Mostowski models for $\Phi$ exist.\) \roster \item"{(i)}" $U$ is a set of
indiscernibles in $\fSg M {K\*U}$ for the  atomic formulas of $\,\Phi$ \(with respect to
the atomic formulas of $\fU$\)\semicolon  in particular\comma if $\bar a$ and $\bar b$ are
two atomically equivalent sequences  in $\fU$\rmp, then they satisfy the same atomic
formulas in $\fSg M {K\*U}$. \item"{(ii)}" $U$ is a set of indiscernibles in $\fSg M
{K\*U}$ for some complete Ehren\-feucht-Mostowski set compatible with $T$\comma 
i.e.\comma if $\bar a$ and $\bar b$ are two atomically equivalent sequences  in $\fU$\rmp,
then they are elementarily equivalent in $\fSg M {K\*U}$. \item"{(iii)}" If $\vartheta$
embeds $\fV$ into $\fU$\comma then the canonical extension of $\vartheta$ to $\Sg N
{(K\*V)}$ elementarily embeds $\fSg N  {K\*V}$ into  $\fSg M {K\*U}$. \endroster
In particular\colon

\roster\item"{(iv)}" Every automorphism of $\fU$ extends canonically to an
automorphism of $\fSg M {K\*U}$.
\item"{(v)}" If $\fV\subseteq\fU$, then $\fSg M  {K\*V}\prc\fSg M
 {K\*U}$. \endroster
The lemma continues to hold if we replace \rmp{``\,}$\fM$\rmp" and
\rmp{``\,}$\fN$\rmp" everywhere by
  \rmp{``\,}$\fM_K$\rmp" and \rmp{``\,}$\fN_K$\rmp" respectively. 
\endproclaim
\demo{Proof} Since satisfaction of atomic formulas is preserved under
submodels, part (i) is 
 trivial.  Part (i) directly implies (iv),
and part (v) implies (iii).  Thus, we only have to prove (ii) and (v). 
We begin with (v).

First, assume that $\fU$ is $\omega$-homogeneous.  We easily check that the
Tarski-Vaught criterion holds.  Here are the details.  Suppose $\bar a$ is an
$n$-termed sequence from $\fSgx {K\*V}$ satisfying ${\exists}x\varphi(x,\bar y)$
in $\fSgx {K\*U}$; say, $b$ is an element of $\fSgx {K\*U}$ such that 
$\langle b \rangle\hat{\ }\bar a$ satisfies $\varphi(x,\bar y)$.  Let $\bar
c$ and $\bar d$  be sequences of elements from $\fV$ and $\fU$ 
that generate $\bar a$ and $b$ respectively.  Thus, there are
$f_0,\dots,f_{k-1}$ and $g_0,\dots,g_{l-1}$ in $K$,  and a term $\sigma(\bar z,
\bar w)$ and terms $\tau_0(\bar u, \bar v),\dots,\tau_{n-1}(\bar u, \bar v)$
in $L$ such that (adding dummy variables to simplify
notation)
$$a_i=\tau_i\bigl[f_0[\bar c],\dots,f_{k-1}[\bar c],\bar c\bigr]\quad\text{for $i
< n$, and}\quad b = \sigma\bigl[g_0[\bar d],\dots,g_{l-1}[\bar d], \bar d\,\bigr] \quad.$$ 
Hence, the sequence 
 $$\bigl\langle\,\sigma\bigl[g_0[\bar d\,],\dots, g_{l-1}[\bar d\,],\bar
d\,\bigr],\tau_0\bigl[f_0[\bar c],\dots,f_{k-1}[\bar c], \bar c\bigr],
\ldots,\tau_{n-1}\bigl[f_0[\bar c],\dots,f_{k-1}[\bar c],\bar
c\bigr]\,\bigr\rangle$$
 satisfies $\varphi$ in $\fSgx
{K\*U}$. Since $\fV$ is dense, there is a sequence  $\bar e$ from $\fV$ such
that $\bar c\,\hat{\ }\bar d$ is atomically equivalent to  $\bar c\,\hat{\ }\bar
e$ in $\fU$.  By $\omega$-homogeneity there is an automorphism  of $\fU$
taking $\bar c\,\hat{\ }\bar d$  to $\bar c\,\hat{\ }\bar e$. In view  of (iv),
this automorphism extends to an automorphism of $\fSgx {K\*U}$.  Thus,
$$\bigl\langle\,\sigma\bigl[g_0[\bar e],\dots, g_{l-1}[\bar e],\bar
e\bigr],\tau_0\bigl[f_0[\bar c],\dots,f_{k-1}[\bar c], \bar c\bigr],
\ldots,\tau_{n-1}\bigl[f_0[\bar c],\dots,f_{k-1}[\bar c],\bar
c\bigr]\,\bigr\rangle$$ satisfies $\varphi$ in $\fSgx {K\*U}$. Since
$\sigma[g_0[\bar e],\dots, g_{l-1}[\bar e],\bar e]$ is in 
$\fSgx {K\*V}$, this shows that the Tarski-Vaught criterion is satisfied. 
Hence,  $\fSgx {K\*V} \prc\fSgx {K\*U}$.

Now let $\fU$ be an arbitrary dense ordering or dense tree.  Take $\fW$ to be an
$\omega$-homogeneous extension of $\fU$, and set $\fP'=\EM W$.  Thus,
$\fM'\prc\fP'$, and hence $\fM\prc\fP$, by Fact 3.  By the case just 
treated we have 
$$\fSg M {K\*U}=\fSg P {K\*U}\prc \fSg P {K\*W}\quad\,$$
and $$\fSg M {K\*V}=\fSg N {K\*V}\prc \fSg N {K\*W}\quad.$$ 
Hence, $\fSg M {K\*V}\prc \fSg M {K\*U}$, as was to be shown.  This proves (v).

For (ii), suppose $\bar a,\bar b\,\in{ }^n U$ are atomically equivalent in
$\fU$.  Let $\fW$ be an $\omega$-homogeneous elementary extension of $\fU$ that
includes the elements of $\bar a$ and $\bar b$, and set $\fP'=\EM W$.  By
Fact 3 we may assume that $\fM'\prc\fP'$, and hence that $\fM\prc\fP$. 
Since $\fW$ is $\omega$-homogeneous,
there is an automorphism of
 $\fW$ taking $\bar a$
to $\bar b$.  By Fact 2, this extends to an automorphism $\vartheta$   of
$\fP'$. Now $\vartheta$ maps $W$ one-one onto itself and preserves all the
operations of $\fP'$.  In particular, an appropriate restriction is an
automorphism of $\fSg P {K\*W}$.  Thus, $\bar a$ and $\bar b$ are elementarily
equivalent in $\fSg P {K\*W}$. By part (v) they remain elementarily equivalent
in $\fSg P {K\*U} = \fSg M {K\*U}$.

The above proof obviously remains valid if we replace
(implicit and explicit occurrences of) \rmp{``\,}$\fM$\rmp" and
\rmp{``\,}$\fN$\rmp" everywhere by
  \rmp{``\,}$\fM_K$\rmp" and \rmp{``\,}$\fN_K$\rmp" respectively.
\qed\enddemo

\redefine\fSg#1#2{\frak {Sg}^{\frak{#1}}#2}
\redefine\fSgx#1{\frak {Sg}@,@,@,#1}

\proclaim{Lemma 2}  Suppose $\fU$ has power $\ka$\comma  $\fM'=\EM U$\comma and
$K\subseteq L'$ has power at most $\ka$.  Then $\fSg M {(K\*U)}$ is strictly
$\ka$-generated.
\endproclaim
\demo{Proof} Since $|U|=\ka\leq |K|$, we have $|K\*U|=\ka$. Thus,
$\fSgx {(K\*U)}$ is   $\ka$-generated.  Suppose now, for contradiction,
that $\fSgx {(K\*U)}$ is $\mu$-generated for some (finite or
infinite) $\mu < \ka$.  A standard argument then gives us a set
$V\subseteq U$, of power $\mu$ when $\omega\leq\mu$, and finite when
$\mu < \omega$, such that $K\*V$ generates $\fSgx {(K\*U)}$. In
particular, $K\*V$ generates $U$ in $\fM$.  But then $V$ generates $U$ in
$\fM'$, which is impossible, by Fact 4.  Indeed, either $U$ is irredundant in
$\fM'$ or else a cardinality argument gives us an $f$ in $U$ that is not an initial
segment of any element of $V$.\qed\enddemo 

\proclaim{Main Hypothesis} Throughout the remainder of the paper we fix an
uncountable cardinal  $\ka$ and assume that $T$  is
a $\ka$-unique universal theory in a language $L$ of arbitrary
cardinality.\endproclaim
As before, $L'$
will be an arbitrary expansion of $L$ with built-in Skolem functions, $T'$ a
Skolem theory in $L'$ extending $T$, and $\Phi$ a standard or generalized
Ehrenfeucht-Mostowski set compatible with $T'$.  Again, recall that
if $\fM' = \EM U $, then we set $\fM =  \rst {\fM'} L$. Our first goal is to show that
$T_\infty$---the theory of the infinite models of $T$---is complete.

\proclaim{Lemma 3} Let $\fU$ be a dense ordering of
power $\ka$\comma and  $\fM'=\EM U$. Then for every $n$ and every countable
set $H\subseteq L'-L$\comma there is countable set $K$\comma with $H\subseteq
K\subseteq L'-L$\comma such that  $$ \fSg M {(K\*U)} \prc_n \fM\quad.$$ The same is
true if $\fU$ is a dense tree\comma provided that $\fM'$ exists.  \endproclaim
\demo{Proof} The proof is by induction on $n$, for all $H$ at once. The case
$n=0$ is trivial: take $K=H$. Assume, now, that the lemma holds for a given
$n\geq 0$ and for all countable sets $H\subseteq L'-L$. Suppose, for
contradiction, that $H_0$ is a countable subset of $L'-L$ such that, 
\num 1 {For every countable $K$ with $H_0\subseteq K\subseteq L'-L$ we
have}
$$\fSg M {(K\*U)} \nprc_{n+1} \fM\quad.$$
We construct an increasing sequence $\langle\,H_\xi :
\xi<\omega_1\,\rangle$  of countable subsets of $L'-L$ such that, setting 
$$\fA_\xi = \fSg M {(H_\xi\*U)}\,,$$
we have
\num 2 {$\fA_\xi \prc_n \fM\quad\text{ and }\quad
\fA_\xi\nprc_{n+1} \fA_{\xi+1}\quad\text{ for
all }\ \xi<\omega_1$\quad.} 
\endnum
The set $H_0$ is given.  Suppose $H_\xi$ has been constructed.  Since $\fA_\xi
\nprc_{n+1}\fM$, by (1), there is a $\Pi_n$-formula
$\varphi_\xi(x_0,\ldots,x_{p-1},y_0,\ldots,y_{q-1})$ in $L$, appropriate terms
$\tau_0,\ldots,\tau_{q-1}$ in $L$ and function symbols $f_0,\ldots,f_{r-1}$ in
$H_\xi$\p, and a sequence $\bar a=  \langle\,a_0,\ldots,a_{s-1}\,\rangle$ 
from $U$ such that (adding dummy variables to simplify notation) 
\medskip \smallskip\flushpar
 \halign{#\hfill&\quad#\hfill\cr
(3)&$\langle\,f_0[\bar a],\ldots,f_{r-1}[\bar a], a_0,\ldots,a_{s-1}\rangle
$ satisfies
$\exists{x_0}\ldots\exists{x_{p-1}}\varphi_\xi(x_0,\ldots,x_{p-1},
\tau_0,\ldots,\tau_{q-1})$ \cr
\ & in $\fM$,  but not in $\fA_\xi$\quad.\cr}  
\medskip \smallskip\flushpar
(Here we are
using the fact that $H_\xi\*U$ generates $\fA_\xi$\p.) Thus, there are Skolem
functions $g_0,\dots,g_{p-1}$ in $L'$ such that 
\medskip \smallskip\flushpar
 \halign{#\hfill&\quad#\hfill\cr
(4)&$\langle\,g_0[\bar a],\ldots,g_{p-1}[\bar a],f_0[\bar a], \ldots,
f_{r-1}[\bar a],a_0,\ldots,a_{s-1}\,\rangle$ satisfies the $\Pi_n$-formula \cr
\ & $\varphi_\xi(x_0,\ldots,x_{p-1},
\tau_0,\ldots,\tau_{q-1})$ in $\fM$\quad.\cr}  
\medskip \smallskip\flushpar
By the induction hypothesis there is a
countable $H_{\xi+1}$ with $$H_\xi\cup\{g_0,\ldots,g_{p-1}\}\subseteq
H_{\xi+1}\subseteq L'-L$$ and such that, setting
$\fA_{\xi+1}=\fSgx {(H_{\xi+1}\*U)}$, we have
$\fA_{\xi+1}\prc_n\fM$.
However, we put witnesses in $\fA_{\xi+1}$ to insure that
$\fA_{\xi}\nprc_{n+1}\fA_{\xi+1}$\p.
At limit stages $\delta$ take $H_\delta=\bigcup_{\xi<\delta}H_\xi$\p. This
completes the construction.

Set $$G={\tsize\bigcup}_{\xi<\omega_1}H_\xi\quad\text{ and }\quad
\fB=\fSg M {(G\*U)}\quad.$$
Clearly, $\fB=\bigcup_{\xi<\omega_1}\fA_\xi$\p.  Therefore, from (2) and our
construction we obtain
\num 5 {$\fA_\xi\prc_n\fB\prc_n\fM\quad\text{ and }\quad 
\fA_\xi\nprc_{n+1}\fB\quad\text{ for every }\xi<\omega_1\quad.$}
\endnum

Next, we prove:
\num 6 {For every $\xi< \omega_1$ and every countable, 
dense $X\subseteq U$
we have}
$$\fSgx {(H_\xi\*X)} \nprc_{n+1}\fB\quad.$$
Indeed, fix $\xi < \omega_1$ and select a sequence $\bar a$ 
so that (3) and (4) hold. Given $X\subseteq U$, choose a sequence $\bar b$ from $X$
that is atomically equivalent to $\bar a$ in $\fU$. This is possible because
$X$ is dense. By indiscernibility in $\fM'$, and (4), $$\langle\,g_0[\bar
b],\ldots,g_{p-1}[\bar b],f_0[\bar b], \ldots,f_{r-1}[\bar
b],b_0,\ldots,b_{s-1}\,\rangle$$ satisfies
$\varphi_\xi(x_0,\ldots,x_{p-1},\tau_0,\ldots,\tau_{q-1})$ in $\fM$, and hence
also in $\fB$, by (5).  Suppose, for contradiction, that $\fSgx {(H_\xi\*X)}
\prc_{n+1}\fB$. Then  \medskip \smallskip\flushpar
 \halign{#\hfill&\quad#\hfill\cr
(7)&$\langle f_0[\bar b],\ldots,f_{r-1}[\bar b],
b_0,\ldots,b_{s-1}\rangle\text{ satisfies }\exists{x_0}\ldots\exists{x_{p-1}}
\varphi_\xi(x_0,\ldots,x_{p-1},\tau_0,\ldots,\tau_{q-1})$\cr
\ & in $\fSgx {(H_\xi\*X)}$\quad.\cr}  
\medskip \smallskip\flushpar
Now we have $\fSgx {(H_\xi\*X)} \prc \fA_\xi$\p, by Lemma 1(v).  Therefore,
(7) holds with ``$\fSgx {(H_\xi\*X)}$" replaced by ``$\fA_\xi$".  Since
$\bar a$ and $\bar b$ are elementarily equivalent in the expansion  $\frak
{Sg}^{\frak{M}_{H_\xi}} {(H_\xi\*U)}$ of $\fA_\xi$\p, by Lemma 1(ii), we
conclude that  $\langle f_0[\bar a],\ldots,f_{r-1}[\bar a],
a_0,\ldots,a_{s-1}\rangle$  satisfies
$$\exists{x_0}\ldots\exists{x_{p-1}}
\varphi_\xi(x_0,\ldots,x_{p-1},\tau_0,\ldots,\tau_{q-1})$$ in $\fA_\xi.$
But this contradicts (3); so the proof of (6) is completed. 

We now work towards a contradiction of (6) by producing a $\xi < \omega_1$ and a
countable, dense set $X\subseteq U$ such that 
\num 8 {$\fSgx {(H_\xi\*X)} \prc \fB\quad.$}
\endnum
Since $\fB=\fSg M {(G\*U)}$ and $\fSg M {(U)}$ are both strictly
$\ka$-generated, by Lemma 2, they are isomorphic, by $\ka$-uniqueness.  Thus,
there is a set $V\subseteq B$ and a dense ordering $\sqsubseteq$ on $V$ such that
$\fV =\lr V \sqsubseteq$ and $\fU$ are isomorphic, $V$ generates
$\fB$, and $V$ is a set of indiscernibles in $\fB$ with respect to
 atomic formulas (under the ordering $\sqsubseteq$).

We define increasing sequences, $\langle\,X_n:n \in \omega\,\rangle$ and 
$\langle\,Y_n:n \in \omega\,\rangle$, of countable, dense subsets of $\fU$ and $\fV$
respectively, and an increasing sequence 
$\langle\,\rho_n:n \in \omega\,\rangle$ of countable ordinals, such that 

\num 9 {$\fSg B {(H_{\rho_n}\*X_n)}\subseteq\fSg B {(Y_n)}\subseteq
\fSg B {(H_{\rho_{n+1}}\*X_{n+1})}\quad.$ }
\endnum
Indeed, set $\rho_0=0$ and take $X_0$ to be an arbitrary countable, dense
subset of $U$. Since $H_0$ is countable, so is $H_0\*X_0$\p.  Therefore,
there is a countable $Y_0\subseteq V$ that generates $H_0\*X_0$ in $\fB$. 
By throwing in extra elements, we may assume that $Y_0$ is dense.

Now suppose that $\rho_n$\p, $X_n$ and $Y_n$ have been defined.  Since $Y_n$ is
countable and \linebreak $\fB=\bigcup_{\xi < \omega_1}\fA_\xi$\p, there is a 
$\rho_{n+1} < \omega_1$ such that $Y_n\subseteq A_{\rho_{n+1}}$\p. 
Without loss of generality we may take
$\rho_{n+1} > \rho_n$\p.  Since 
$H_{\rho_{n+1}}\*U$ generates
$\fA_{\rho_{n+1}}$\p, there is  a countable subset
$X_{n+1}\subseteq U$ such that
$H_{\rho_{n+1}}\*X_{n+1}$ generates  $Y_n$ in
$\fB$.  Of course we may choose $X_{n+1}$ so that it is dense and includes
$X_n$\p. As before, we now choose a
countable, dense $Y_{n+1} \subseteq V$ that generates
$X_{n+1}$ in $\fB$ and includes $Y_n$\p.  This completes
the construction.

Set
$$\rho=\sup\{\rho_n:n\in\omega\}\quad, \quad X={\tsize\bigcup}_{n\in\omega}X_n\quad,
\quad Y={\tsize\bigcup}_{n\in\omega}Y_n\quad.$$
Then $\rho<\omega_1$ and $H_\rho=\bigcup_{n\in\omega}H_{\rho_n}$\p, by
definition of the sequence $\langle\,H_\xi:\xi<\omega_1\,\rangle$.  Also, $\lr X \leq$
 and $\lr Y \sqsubseteq$ are countable, dense submodels of $\fU$ and $\fV$
respectively.  By (5), (9), and Lemma 1(v) we have 
$$\fSg M {(H_\rho\*X)}=\fSg B
{(H_\rho\*X)}={\tsize\bigcup}_{n <\omega} \fSg B
{(H_{\rho_n}\*X_n)}={\tsize\bigcup}_{n<\omega} \fSg B
{(Y_n)}=\fSg B {(Y)} $$ and
$$\fSg B {(Y)}=\fSg M {(Y)}\prc\fSg M {(V)}=\fSg B {(V)}=\fB\quad.$$
Thus, we have constructed a countable $\xi<\omega_1$ and a countable, dense
$X\subseteq U$ such that (8) holds.  This contradicts (6).
\qed\enddemo
\proclaim{Lemma 4} 
For
every countable $H\subseteq L'-L$\comma there is a countable $K$\comma with
$H\subseteq K\subseteq L'\diff L$\comma such that\comma whenever $\fV$ is a dense
ordering and $\fN'=\EM V$\comma we have $$\fSg N {(K\*V)} \prc \fN\quad.$$  The same
is true when $\fV$ is a dense tree\comma provided that $\EM U$ exists for dense trees
$\fU$. \endproclaim  \demo{Proof} Let $\fU$ be a dense ordering of power $\ka$, and set
$\fM'=\EM U$.  We use
Lemma 3 to define an increasing sequence $\langle\,H_n :n<\omega\,\rangle$
of countable subsets of $L'-L$ such that \num 1 {$\fSg M {(H_n\*U)} \prc_n
\fM\quad.$} \endnum
Indeed, we set $H_0=H$, and given $H_n$ satisfying (1), we apply Lemma 3 to
obtain $H_{n+1}$\p. Setting $K=\bigcup_{n\in\omega}H_n$\p, we easily check 
that 
\num 2 {$\fSg M (K\*U)\prc\fM$\quad.}
\endnum

Now suppose that $\fV$ is any dense ordering or tree, and that $\fN'=\EM
V$.  We first treat the case when $\fV$ is a submodel of $\fU$.  By Fact 3
we may assume, without loss of generality, that $\fN'\prc\fM'$.  Hence, 
\num 3 {$\fN\prc\fM$\quad.}
\endnum
In view of (3) and Lemma 1(v), we get
\num 4 {$\fSg N {(K\*V)} = \fSg M {(K\*V)}\prc \fSg M {(K\*U)}$\quad.}
\endnum
Combining (2)--(4) gives
\num 5 {$\fSg N {(K\*V)}\prc \fN$\quad.}
\endnum

For the case when $\fV$ is not a submodel of $\fU$, take a countable, dense submodel
$\fW$ of $\fV$.  Set $\fP'=\EM W$.  Again, we may assume that $\fP'\prc\fN'$, so
\num 6 {$\fP\prc\fN$\quad,}
\midnum 7 {$\fSg P {(K\*W)} = \fSg N {(K\*W)}\prc \fSg N {(K\*V)}$\quad.}
\endnum
Furthermore, there is an embedding of $\fW$ into $\fU$, and this embedding
extends to an elementary embedding $\vartheta$ of $\fP'$ into $\fM'$.  Hence,
as usual,
\num 8 {$\vartheta$ induces an elementary embedding of $\fP$ into $\fM$\quad,}
\midnum 9 {$\vartheta$ induces an elementary embedding of $\fSg P {(K\*W)}$
into  $\fSg M {(K\*U)}$\quad.}
\endnum

Using (2) and (6)--(9), we readily verify that (5) holds.  Here are the details. 
Let $\varphi(\bar x)$ be a formula of $L$, and $\bar a$ a sequence of appropriate
length from $\fSg N {(K\*V)}$.  Then there is a finite sequence $\bar v$ from
$V$ that generates $\bar a$ (with the help of $K$).  Choose a $\bar v'$ in $W$
that is atomically equivalent to $\bar v$, and let $\bar a'$ be the sequence
obtained from $\bar v'$ in the same way that $\bar a$ is obtained from $\bar
v$.  Then $\bar a$ and $\bar a'$ are elementarily equivalent in $\fN'$---and
therefore also in $\fN$---and in $\fSg N {(K\*V)}$ by Fact 3 and Lemma 1(ii). 
Hence,  $$\alignat 2
\fN\models\varphi[\bar a]\qquad&\text{iff}\qquad \fN\models\varphi[\bar a']&&\\
&\text{iff}\qquad \fP\models\varphi[\bar a']&&\qquad\text{by (6),}\\
&\text{iff}\qquad \fM\models\varphi[\vartheta(\bar a')]&&\qquad\text{by (8),}\\
&\text{iff}\qquad \fSg M {(K\*U)}\models\varphi[\vartheta(\bar
a')]&&\qquad\text{by (2),}\\ 
&\text{iff}\qquad \fSg P {(K\*W)}\models\varphi[\bar
a']&&\qquad\text{by (9),}\\
&\text{iff}\qquad \fSg N {(K\*V)}\models\varphi[\bar
a']&&\qquad\text{by (7),}\\ 
&\text{iff}\qquad \fSg N {(K\*V)}\models\varphi[\bar
a]\quad.&&
\endalignat
$$
 \qed\enddemo

\proclaim{Theorem 5} $T_\infty$ is complete.
\endproclaim
\demo{Proof} Let $\fB_1$ and $\fB_2$ be two infinite models of $T$. For 
$i=1,2$, let $T_i$ be the  theory of $\fB_i$\p, let $L'$ be an expansion of the
language of $T_i$ with built-in Skolem functions, $T'_i$ a Skolem theory in $L'$
extending $T_i$\p, and $\Phi_i$ a standard Ehrenfeucht-Mostowski set in $L'$
compatible with $T'_i$\p. Let $\fU$ be a dense ordering of power $\ka$, and set
$\fM'_i=EM(\fU,\Phi_i)$ and $\fM_i = \rst{\fM'_i}{L}$. By Lemma 4 there is
a countable $K_i\subseteq L'-L$ such that 
\num 1 {$\frak {Sg}^{\fM_i}(K_i\*
U)\prc\fM_i$\quad.}
\endnum
Now $\frak {Sg}^{\fM_1}(K_1\*
U)$ and  $\frak {Sg}^{\fM_2} (K_2\*U)$ are strictly $\ka$-generated
models of $T$, by Lemma 2.  Therefore, they are isomorphic, by the
$\ka$-uniqueness of $T$.  From this and (1) we conclude that $\fM_1$
and  $\fM_2$ are elementarily equivalent.  But $\fM'_i$\,, and hence also  
$\fM_i$\,, is a
model of the theory of $\fB_i$\p.  In consequence, $\fB_1$ and $\fB_2$ are
elementarily equivalent, as was to be shown.
\qed\enddemo

Our next goal is to prove:
\proclaim{Theorem 6} $T_\infty$ is superstable. \endproclaim
\demo{Proof}  Set $U=(\kappa + 1)\times\Bbb Q$ and $\fU = \lr U <$, where $<$
is the lexicographic ordering on $U$. We define a substructure $\fV$ of the
full tree over $U$ by specifying its universe. It
is the smallest set $V$ satisfying the following conditions: all finite levels
of the full tree are included in $V$, i.e., $\mybigcup_{n\in \omega}\!{}^n U
\subseteq V$; all eventually constant functions from ${}^\omega U$ are in $V$;
for each limit ordinal $\delta<\omega_1$ we choose a strictly increasing function
$f_\delta$ in ${}^\omega{\omega_1}$ such that
$\sup\{f_\delta(n):n\in\omega\}=\delta$, and we put into $V$ a copy of $f_\delta$ from 
${}^\omega V$ called
$g_\delta$ and determined by: $g_\delta(n) =
\lr{f_\delta(n)} 0$ for each $n$.
\num 1 {For every countable set
$W\subseteq U$, the set $V\cap\preom W$ is also countable.}
\endnum

Indeed, the finite levels, ${}^n W$, and the set of eventually constant
functions of ${}^\omega W$ are all countable.  Moreover, the domain of $W$,
i.e., $\{\delta: \lr\delta q \in W \text{ for some } q\in \Bbb Q\}$, is
countable, so the set of $g_\delta$ in ${}^\omega W$, with $\delta<\omega_1$,
is also countable.  This proves (1).

Let $\Phi$ be a standard Ehrenfeucht-Mostowski set compatible with $T'$ (our Skolem
extension of $T$\/), and set $\fM' = \EM U$.  Assume, for contradiction, that
$T$ is not superstable, and let $\langle \varphi_n(x,y):n<\omega\rangle$ be a
sequence of formulas witnessing this nonsuperstability.  Then, as mentioned in
the preliminaries, there is an Ehrenfeucht-Mostowski set $\Psi$ compatible with
$T'$ such that the generalized Ehrenfeucht-Mostowski model $\fN'=\EMM V$ 
exists, and, for $u$ from $P_n^\fV$ and $v$ from $P_\omega^\fV$, we have
$\fN'\models \varphi_n[u,v]$ iff $u\lessdot v$.  Set $\fM = \rst {\fM'} L$ and 
$\fN = \rst {\fN'} L$.

By Lemma 4, there are countable sets $J,K\subseteq L'\diff L$ such that
\num 2{$\fSgx {(J\* U)}\prc \fM$\quad and \quad
$\fSgx {(K\* V)}\prc \fN$\quad.}
\endnum
By Lemma 2, $\fSgx {(J\* U)}$ and
$\fSgx {(K\* V)}$ are both strictly $\kappa$-generated.  Hence, $\kappa$-uniqueness
implies that they are
isomorphic.  Let $\vartheta$ be such an isomorphism.

We now define two strictly increasing, continuous sequences, $\langle
X_\alpha:\alpha <\omega_1\rangle$ and \linebreak
$\langle
Y_\alpha:\alpha <\omega_1\rangle$, of countable, dense subsets of $U$ and $V$
respectively,  each $Y_\alpha$  being downward closed, i.e., closed under
initial segments.  (i)  Let $X_0 = \{\kappa\}\times\Bbb Q$ and take $Y_0$ to be
any countable, dense, downward closed subset of $V$.  Suppose, now, that
$X_\alpha$ and $Y_\alpha$ have been defined.  Since both of these sets are
countable, by assumption, so are $\jx$ and $\ky$.  For any $a$ in 
$J\*U$ there is obviously a finite subset $C_a\subseteq V$ such that $K\*C_a$
generates $\vartheta(a)$ (in $\fSgx (K\*V)$).  (ii) For each $a$ in
$\jx$, put the elements of  $C_a$ into $Y_{\alpha + 1}$.  Since $\jx$ is
countable, this adds only countably many elements to $Y_{\alpha + 1}$. 
(iii)  Similarly, for every $b$ in $\ky$, choose a finite set $D_b\subseteq U$
such that $J\* D_b$ generates $\vartheta^{\sssize -1}(b)$ (in $\fSgx (J\*U)$),
and  put each element of $D_b$ into $X_{\alpha + 1}$.
(iv) For each $\lr \xi q$ in $X_\alpha$, put all of $\{\xi\}\times\Bbb Q$ into
$X_{\alpha + 1}$.  (v)  Put all of $\preom (\alpha\times\Bbb Q)\cap V$ into
$Y_{\alpha + 1}$. The latter set remains countable, by (1).  (vi) If necessary,
add countably many more elements to $X_{\alpha + 1}$ and $Y_{\alpha + 1}$ to
insure that these sets are dense, that $X_\alpha\subset X_{\alpha + 1}$ and
$Y_\alpha\subset Y_{\alpha + 1}$, and that $Y_{\alpha +1}$ is downward closed.
(vii)  For limit ordinals $\delta$, set $X_\delta = \mybigcup_{\alpha <\delta}
X_\alpha$ and $Y_\delta = \mybigcup_{\alpha <\delta}
Y_\alpha$. 

Setting 
\num 3 {$\fA_\alpha =\fSgx
(\jx)$ and  $\fB_\alpha =\fSgx (\ky)$ for every $\alpha < \omega_1$\quad,}
\endnum
we see from conditions (ii) and (iii) that, for each  $\alpha<\omega_1$, we
have $$\vartheta(A_\alpha)\subseteq B_{\alpha +1}\quad \text{and} \quad
\vartheta^{\sssize -1}(B_\alpha)\subseteq A_{\alpha +1}\quad.$$
\relax From this it easily follows that
\medskip\smallskip\flushpar 
\halign{#\hfill&\quad#\hfill\cr
(4)&For each limit ordinal $\delta <\omega_1$, the (appropriate restriction
of the) mapping $\vartheta$ is an\cr
\ & isomorphism of $\fA_\delta$ onto
$\fB_\delta$\quad.\cr}  
\medskip \smallskip
Set $$E=\bigl\{\delta<\omega_1:\text{$\delta$ is a limit ordinal and }\preom
(\omega_1\times\Bbb Q) \cap Y_\delta = \bigl[\mybigcup_{\beta<\delta}
\preom(\beta\times\Bbb Q)\bigr]\cap V\bigr\}$$
and
$$E'=\bigl\{\delta\in E: \delta =\sup (\delta\cap E) = \sup\{\alpha\in E: \alpha
<\delta\}\bigr\}\quad.$$ 
\num 5 {$E$ and $E'$ are closed, unbounded sets.}
\endnum

To see that $E$ is unbounded, observe, first of all, that, by conditions (v)
and (vii), we have 
\num 6 {$\preom
(\omega_1\times\Bbb Q) \cap Y_\delta \supseteq \bigl[\mybigcup_{\beta<\delta}
\preom(\beta\times\Bbb Q)\bigr]\cap V$\quad for every limit ordinal
$\delta<\omega_1$\quad.} \endnum
Therefore, we need only establish the reverse inclusion for unboundedly many
limit ordinals.  Since a $v$ in $\preom (\omega_1\times\Bbb
Q)\cap V$ must have a countable range, there is an ordinal $\gamma_v<\omega_1$
such that $v$ is in $\preom (\gamma_v\times\Bbb Q)$.  Let $\alpha_0$ be an
arbitrary countable ordinal.  Since $Y_{\alpha_0}$ is countable,
the supremum of the $\gamma_v$, over all $v$ in $Y_{\alpha_0}$, is
countable. Hence, we can find a countable $\alpha_1 > \alpha_0$ such that
$$\preom (\omega_1\times\Bbb Q)\cap Y_{\alpha_0}\subseteq \preom
(\alpha_1\times\Bbb Q)\quad.$$ 
Continuing in this fashion, we obtain a strictly increasing sequence $\langle
\alpha_n:n < \omega\rangle$ of count\-able ordinals such that 
$$\preom
(\omega_1\times\Bbb Q)\cap Y_{\alpha_n}\subseteq \preom
(\alpha_{n+1}\times\Bbb Q)\quad.$$ 
Set $\delta = \sup\{\alpha_n:n <\omega\}$.  Then 
$$\preom
(\omega_1\times\Bbb Q)\cap Y_{\delta}\subseteq
\bigl[\mybigcup_{\beta<\delta}\preom (\beta\times\Bbb Q)\bigr]\cap V\quad.$$ 
In view of (6), we see that equality actually holds in
the preceding line.  This shows that $E$ is unbounded.  The rest of the proof
of (5) is easy, so we leave it to the reader.

\num 7 {For every $\delta$ in $E'$ and every $\alpha$ in $E\cap\delta$, the type
$tp(g_\delta,B_\delta)$ splits over $B_\alpha$\quad.}
\endnum

To prove (7), recall that $g_\delta$ is not in 
${}^\omega (\beta\times\Bbb Q)$\ for any $\beta <\delta$.  We thus see, from the
definition of $E$, that $g_\delta$ is not in $Y_\delta$.  Hence, by Fact 4,
$g_\delta$ is not in $B_\delta$.

The function $g_\delta $ is defined in terms of a strictly increasing,
ordinal-valued function $f_\delta$ whose supremum is $\delta$.  Since
$\alpha <\delta$, there must be an $n <\omega$ such that $f_\delta(n) >
\alpha$.  Let $h_1$ and $h_2$ be the extensions of $\rst {g_\delta} n$ to $n+1$
determined by  $$h_1(n) = g_\delta(n) = \lr {f_\delta(n)}
0\qquad\text{and}\qquad h_2(n) = \lr {f_\delta(n)} 1\quad.$$ 
Notice that $h_1 = \rst {g_\delta} {(n+1)}$.  Clearly, $h_1$ and $h_2$ are in
$Y_\delta$, by conditions (v) and (vii).  Because $\alpha$ is in $E$, and
$f_\delta(n) >\alpha$, we can apply the definition of $E$ to conclude that
$h_1$ and $h_2$ are not in $Y_\alpha$, and hence are not initial segments of
any elements in $Y_\alpha$.  Therefore, they are not in $B_\alpha$, by Fact
4.  They realize the same type over $B_\alpha$, by tree indiscernibility. 
Moreover,  $\fN\models \varphi_{n+1}[h_1,g_\delta]$ and 
$\fN\models \neg\varphi_{n+1}[h_2,g_\delta]$, since $\varphi_{n+1}$ codes
the initial segment relation between elements in $P^{\fV}_{n+1}$ and 
$P^{\fV}_{\omega}$.  Thus, both 
$\varphi_{n+1}(h_1,\bar x)$ and 
$\neg\varphi_{n+1}(h_2,\bar x)$ are in $tp(g_\delta,B_\delta)$.  This
completes the proof of (7).

We now work towards a contradiction to (7).  
\medskip\smallskip
 \halign{#\hfill&\quad#\hfill\cr
(8)&For every $\delta$ in $E'$ and every $a$ in $\Sgx (J\*U)\diff
A_\delta$, there is an $\alpha$ in $E\cap\delta$ such that \cr
\ &  $tp(a,A_\delta)$ does not 
split over $A_\alpha$\quad.\cr}  
\medskip \smallskip\flushpar

To see that (8) contradicts (7), take any $\delta$ in $E'$, and set
$a=\vartheta^{\sssize -1}(g_\delta)$.  Since $g_\delta$ is in 
 $\Sgx (K\*V)\diff
B_\delta$, we get that $a$ is in $\Sgx (J\*U)\diff
A_\delta$, by (4).  Hence, there is an $\alpha$ in $E\cap\delta$ such that
$tp(a,A_\delta)$ does not split over $A_\alpha$, by (8).  But then, applying
$\vartheta$ to $a$, $tp(a,A_\delta)$, and $A_\delta$, we get that
$tp(g_\delta,B_\delta)$ does not split over $B_\alpha$, by (4).  This
contradicts (7).

To prove (8), fix a $\delta$ in $E'$ and an $a$ in $\Sgx (J\*U)\diff
A_\delta$.  Then there is a sequence $\bar u = \langle
u_0,\dots,u_{r-1}\rangle$ from $U$, function symbols $f_0,\dots f_{n-1}$ in
$J$, and a term $\sigma$ from $L$, such that (adding dummy variables)
\num 9 {$a= \sigma\bigl [f_0[\bar u],\dots,f_{n-1}[\bar u],\bar u\bigr]$\quad.}
\endnum

Each $u_i$ has the form $u_i = \lr{\gamma_i} {q_i}$ for a unique $\gamma_i
\leq \kappa$ and $q_i$ in $\Bbb Q$.  The set 
$$\Gamma_i =\bigl\{\beta\leq\kappa:\{\beta\}\times\Bbb
Q\subseteq X_\delta \text{ and } \gamma_i\leq\beta\bigr\}$$
is nonempty, since it contains $\kappa$, by condition (i).  We shall denote
the smallest element of this set by $\gamma_i'$.  Notice that $u_i\in
X_\delta$ iff $\gamma_i' = \gamma_i$, by condition (iv).  Now, by definition,
$\{\gamma_i'\}\times\Bbb Q\subseteq X_\delta$ for each $i<r$.    Since
$\delta$ is in $E'$, there is, for each $i<r$, a $\beta$ in $E\cap\delta$
such that $\{\gamma_i'\}\times\Bbb Q \subseteq X_\beta$.  Take $\alpha$ to be
the maximum of these $\beta$.    Then $\alpha$ is in $E\cap\delta$ and
$\{\gamma_i'\}\times\Bbb Q\subseteq X_\alpha$ for each $i<r$. 

To verify that $tp(a,A_\delta)$ does not split over $A_\alpha$, let
$\psi(x,\bar y)$ be a formula of $L$, and let $\bar b$ and $\bar c$ be
sequences of equal length from $A_\delta$, such that $\psi(x,\bar b)$ and
$\neg\psi(x,\bar c)$ are in $tp(a,A_\delta)$, i.e., 
\num {10} {$\fM\models \psi[a,\bar b]\wedge\neg\psi[a,\bar c]$\quad.}
\endnum
We must prove that 
\num {11} {$tp(\bar b,A_\alpha)\neq tp(\bar c,A_\alpha)$\quad.}
\endnum
Since $\bar b$ and $\bar c$ come from $A_\delta$, there is a finite sequence
$\bar v$ from $X_\delta$ that generates $\bar b$ and $\bar c$ with the help
of some elements of $J$.  Using the ordinals $\gamma_i'$\,, we shall construct a
sequence $\bar u' =\langle u_0',\dots,u_{r-1}'\rangle$ in $X_\alpha$ such
that, for each $i$, we have  $u_i' = u_i$ iff $u_i$ is in $X_\delta$, and
\medskip\smallskip
{\baselineskip=18pt \halign{#\hfill&\quad#&\qquad#&\qquad#&\quad#\hfill\cr
(12)&$u_i'<u_k'$ iff $u_i<u_k$&and &$u_i'=u_k'$ iff $u_i=u_k$&, \cr
\ &$u_i'<v_j$ iff $u_i <
v_j$&and&$u_i'=v_j$ iff $u_i = v_j$&.\cr}}  
\medskip \smallskip\flushpar
for all appropriate $j,k$. 

The construction of $\bar u'$ is not difficult:  we shall set $u_i' =\lr
{\gamma_i'} {q_i'}$, where $q_i'$ is chosen from $\Bbb Q$.  Since
$\{\gamma_i'\}\times \Bbb Q\subseteq X_\alpha$, this assures that $u_i'$ will
be in $X_\alpha$.  If $\gamma_i$ is in $\Gamma_i$, then $\gamma_i'=\gamma_i$.  In
this case, take $q_i'=q_i$, so that $u_i'=u_i$.  Suppose, now, that $v_j =
\lr{\beta_j} {r_j}$.  If $u_i < v_j$, then we have either $\gamma_i<\beta_j$, or
else $\gamma_i = \beta_j$ and $q_i < r_j$.  In the first case, observe that
$\beta_j\in\Gamma_i$\,, since $v_j$ is in $X_\delta$ (here we are using again
condition (iv)).  Therefore, $\gamma_i'\leq \beta_j$.  If $\gamma_i'<\beta_j$,
then we may choose $q_i'$ however we wish.  If $\gamma_i'=\beta_j$, then we must
choose $q_i' <r_j$.  This is possible since $\{\gamma_i'\}\times\Bbb Q \subseteq
X_\delta$.  In the second case, when $\gamma_i=\beta_j$, we have $\gamma_i$ in
$X_\delta$, so $u_i'=u_i$.  In both cases we obtain $u_i'<v_j$.  The other cases
in (12) are treated similarly.  Notice, however, that we must choose $q_i'$ so
that (12) holds for all appropriate $j$ and $k$ at once.  To do this, we use the
density of $\Bbb Q$.

\relax From (12) we conclude that $\bar u'\hat{\ }\bar v$ and $\bar
u\hat{\ }\bar v$ 
are atomically equivalent in $\fU$.  Because $U$ is a set of order
indiscernibles in $\fM'$, it follows that 
$\bar u'\hat{\ }\bar v$ and $\bar u\hat{\ }\bar v$ are elementarily equivalent
in $\fM'$. Recalling (9), set 
$a'= \sigma\bigl [f_0[\bar u'],\dots,f_{n-1}[\bar u'],\bar u'\bigr]$, and
notice that $a'$ is in $A_\alpha$, since $\bar u'$ is from $X_\alpha$.

Since $\bar u'$ generates $a'$ in the same way that $\bar u$ generates $a$, and
since $\bar v$ generates $\bar b$ and $\bar c$, we conclude that
$a'\hat{\ }\bar b\hat{\ }\bar c$ and  $a\hat{\ }\bar b\hat{\ }\bar c$ are
elementarily equivalent in $\fM'$.  In particular, from (10) we get
that $$\fM\models\psi[a',\bar b]\wedge\neg\psi[a',\bar c]\quad.$$
Thus, $\psi(a',\bar y)$ is in $tp(\bar b,A_\alpha)$, while 
$\neg\psi(a',\bar y)$ is in $tp(\bar c,A_\alpha)$.  This proves (11), and hence
(8). 
\qed\enddemo

As is well-known, the fact that $T_\infty$ is superstable (and hence
stable) implies that  order indiscernibles are totally indiscernible. Thus,
if $\fM'=\EM U$,  then any two one-one sequences in ${}^n U$ satisfy the same
formulas in  $\fM$ (but not necessarily in $\fM'$).  It follows that Lemma 1
(with $K=\varnothing$) goes through without any references  to order. Thus,
if $\fM'$ and $\fN'$ are as in Lemma 1, then any injection  $\vartheta$ of
$V$ into $U$ extends canonically to an elementary embedding  of $\fSg N
{(V)}$ into $\fSg M {(U)}$, and this extension is an isomorphism  iff
$\vartheta$ is onto.  In particular, the isomorphism type of a subalgebra
(of some very large model) generated by an infinite set of indiscernibles
with respect to $\Phi$  is uniquely determined by the cardinality of the set
of indiscernibles.  To give a precise definition of such algebras, we fix a
very large dense ordering $\fZ$ (as large as we will need for any argument
in this paper), set $\fC=EM(\fZ,\Phi)$, and let $\langle Z_\lambda:
\lambda\text{ an infinite cardinal }\leq |Z|\rangle$ be an increasing
sequence of subsets of $Z$ such that $|Z_\lambda| =\lambda$.
\proclaim{Definition 7}  Let $\fF_\lambda = \fSg C {(Z_\lambda)}$.  
 The set $Z_\lambda$ is referred to as a $\Phi$-basis of
$\fF_\lambda$ \endproclaim  

The algebras $\fF_{\lambda}$\p, for $\lambda\ge\omega$, have many 
properties in common with free algebras. 
For example, if $U$ and $V$ are $\Phi$-bases of $\fF_{\lambda}$
and $\fF_{\mu}$ respectively, and if $\vartheta$ maps $V$ one-one into $U$
then the canonical extension of $\vartheta$ maps $\fF_{\mu}$ elementarily
into $\fF_{\lambda}$\p. This is just a reformulation of what was said above.
For another example, notice that $\fF_{\lambda}$ is strictly \linebreak
$\lambda$-generated, by Lemma 2. It follows that {\it any\/} infinite 
$\Phi$-basis of $\fF_{\lambda}$ must have cardinality $\lambda$; otherwise
we would have $\fF_{\mu}\cong\fF_{\lambda}$ for some $\mu\ne\lambda$, 
forcing the strictly $\mu$-generated model $\fF_{\mu}$ to be 
$\lambda$-generated and visa-versa.

For the case when $T_\infty$ is superstable, we readily extend Lemma 
1(v) to cover $\prc_a:$

\proclaim{Lemma 8} Let $\fU, \fV$ be dense orderings\comma  $\fM' = \EM
U$\comma and $K\subseteq L'\diff L$. If  $\fV\subseteq\fU$\comma then $$\fSg M
{(K\*V)}\prc_a\fSg M {(K\*U)}\quad.$$ 
The lemma continues to hold if we replace \rmp{``\,}$\fM$\rmp" by
  \rmp{``\,}$\fM_K$\rmp". \endproclaim  
\demo{Proof} We begin with the case when $\fU$ is
$\omega_1$-homogeneous.  Let $p$  be a strong type of $\fM$ based on a finite
subset $E$ of $\Sgx {(K\*V)}$, and suppose  that $\pe$ is realized in $\fSgx 
{(K\* U)}$ by $\bar b$. Let $\bar v = \langle v_0,\dots,v_{m-1}\rangle$ be a
finite sequence in $\fV$ that generates $E$ in $\fM$, with the help of finitely
many $f\in K$, and let $\bar u =\langle u_0,\dots,u_{n-1}\rangle$ be a finite
sequence from $\fU$ that generates $\bar b$.  Because of the  density of $\fV$,
we can certainly find a sequence $\bar u'=\langle u'_0,\dots,u'_{n-1}\rangle$ in
$\fV$ such that $\bar v\hat{\ }\bar u$ and $\bar v\hat{\ }\bar u'$ are
atomically equivalent in $\fU$.   Using again the density of $\fV$, is not
difficult to construct an $\omega$-sequence $\bar w$ in $V$ such that  $\bar
v\hat{\ }\bar w\hat{\ }\bar u$ and $\bar v\hat{\ }\bar w\hat{\ }\bar u'$ are
atomically equivalent and the range, $X$, of  $\bar v\hat{\ }\bar w$ is dense. 
Assume, for the moment, that this done. Since $\fU$ is $\omega_1$-homogeneous,
there is an automorphism of $\fU$ taking  $\bar v\hat{\ }\bar w\hat{\ }\bar u$
to $\bar v\hat{\ }\bar w\hat{\ }\bar u'$.  Extend it to an automorphism
$\vartheta$ of $\fM'$.  Then $\vartheta$   also induces an automorphism of
$\fSgx  {(K\* U)}$, and it is the identity mapping on  $\Sgx {(K\*X)}$. Because
$X$ is dense, $\fSgx {(K\*X)}$ is an elementary submodel of $\fSgx {(K\*V)}$, and
therefore $p$ is stationary over this submodel. Since $\vartheta$ fixes an
elementary submodel over which $p$ is stationary, we see that $\vartheta$ must 
map $\pe$ to itself.  Thus, $\vartheta(\bar b)$ realizes $\pe$ in $\fSgx  {(K\*
U)}$. But $\fSgx  {(K\*V)}$ is an elementary substructure of $\fSgx  {(K\* U)}$,
by Lemma 1(v), and all the formulas of $\pe$ are over  $\fSgx  {(K\*V)}$.  Thus,
$\vartheta(\bar b)$ must realize $\pe$ in  $\fSgx  {(K\*V)}$.

Now suppose that $\fU$ is arbitrary.  Let $\fW$ be an $\omega_1$-homogeneous
 extension of $\fU$, and set $\fN' = \EM W$.  Then by the previous
case, 
 $$\fSg M {(K\*V)}=\fSg N {(K\*V)} \prc_a \fSg N {(K\*W)}\quad\,$$
and
 $$\fSg M {(K\*U)}=\fSg N {(K\*U)} \prc_a \fSg N {(K\*W)}\quad .$$
Therefore, $\fSg M {(K\*V)} \prc_a \fSg M {(K\*U)}$.

Returning to $\bar w$, a simple example should suffice to
indicate how it is to be constructed.  Suppose, for example, that 
$$ v_i <  u_0 < u_1 < u_2 < u_3 < v_j\quad,$$ and
that no other elements of the sequences $\bar v$ and $\bar u$ are in the
interval between $ v_i$ and $ v_j$.  Of course, 
$$ v_i <  u'_0 < u'_1 < u'_2 < u'_3 < v_j\quad.$$  If 
$ u'_0 \leq  u_0$, or else if $ u'_3 \geq  u_3$\p, then include
in the range of $\bar w$ a dense set between $ v_i$ and $ u'_0$\p,  or
else between $ u'_3$ and $ v_j$\p, respectively.  Now consider the case when $
u_0 <  u'_0 <  u'_3 <  u_3$.   If 
$ u'_1 \leq  u_1$\p, or else if $ u'_2 \geq  u_2$ then include in
the range of $\bar w$ a dense set between $ u'_0$ and $ u'_1$\p,  or
else between $ u'_2$ and $ u'_3$\p, respectively. There remains the case when $
u_1 <  u'_1 <  u'_2 <  u_2$\p.   In this case,  include in
the range of $\bar w$ a dense set between $ u'_1$ and $
u'_2$.
\qed\enddemo 

The next lemma is the analogue of Lemma 3 for the notion
$\prc_a$.

\proclaim{Lemma 9} Let $\fU$ be a dense
ordering of power $\ka$\comma and $\fM'=\EM U$. Then for every countable $H\subseteq
L'-L$  there is a countable $K$ with $H\subseteq K\subseteq L'-L$ such that
$$\fSg M {(K\*U)}\prc_a \fM\quad.$$ \endproclaim 
\demo{Proof} The proof is quite
similar to the proof of Lemma 3.   Suppose, for contradiction, that $H$ is a
counterexample to the assertion  of the lemma.  Analogously to the proof of
Lemma 3, we construct an increasing, continuous sequence $\langle\,H_\xi :
\xi<\omega_1\,\rangle$ of countable subsets of $L'-L$ such that
\num 1 {$\fSgx  {(H_\xi\*U)}\prc\fSgx 
{(H_{\xi+1}\*U)}\prc\fM$\,,\quad but \quad $\fSgx 
{(H_\xi\*U)}\nprc_a\fSgx  {(H_{\xi+1}\*U)}$\quad.} 
\endnum
Indeed, by Lemma 4 we can
take $H_0$ to be a countable extension of $H$  such that \linebreak $\fSgx
 {(H_0\*U)}\prc\fM$.  Suppose, now,  that $H_\xi$  has been defined, and
that $\fSgx  {(H_\xi\*U)}\prc\fM$.  By the  assumption on $H$,
$\fSgx {(H_\xi\*U)}$ is not $a$-saturated in $\fM$.  Thus, there  is a strong
type $p_{_\xi}$ based on a finite subset $E_\xi$ of $\Sgx {(H_\xi\*U)}$\p,  such
that  $\rst{p_{_\xi}}{^*E_\xi}$ is realized in $\fM$---say by $\bar
a_\xi$---but $\rst{p_{_\xi}}{^*E_\xi}$ is not realized in $\fSgx 
{(H_\xi\*U)}$.  Since $\bar a_\xi$  is generated by $U$ in $\fM'$, there is a
finite set $F\subseteq L'-L$  such that $F\*U$ generates $\bar a_xi$ in
$\fM$. Take $H_{\xi+1}$ to be a  countable extension of $H_\xi\cup F$ in
$L'-L$ such that  $\fSgx 
{(H_{\xi+1}\*U)}\prc\fM$. 
Then $\fSgx  {(H_\xi\*U)}\prc\fSgx  {(H_{\xi + 1}\*U)}$\p.  By construction,
$\bar a_\xi$ realizes $\rst{p_{_\xi}}{{}^*E_\xi}$ in $\fSgx 
{(H_{\xi+1}\*U)}$, so 
$\fSgx  {(H_\xi\*U)}\nprc_a\fSgx  {(H_{\xi + 1}\*U)}$. This completes the
verification of (1).

Set $G=\bigcup_{\xi<\omega_1}H_\xi$ and $\fB=\bigcup_{\xi<\omega_1}
\fSgx  {(H_\xi\*U)}
=\fSgx  {(G\*U)}$.

\num 2 {For every $\xi<\omega_1$  and every countably infinite $X\subseteq U
$ we have $\fSgx  {(H_\xi\*X)}\nprc_a\fB\quad.$}
\endnum

Indeed, let $W$ be a countable, dense subset of $U$ containing  $X$ and
containing finite subsets $U_0$ and $U_1$ of $U$ such that $H_\xi\*U_0$
generates $E_\xi$ and $G\*U_1$ generates $\bar a_\xi$ in $\fM$.  Then, by Lemma
1(v),   \num 3 {$\fSgx  {(G\*W)}\prc \fSgx  {(G\*U)} = \fB$\quad.}
\endnum
Therefore, to prove (2) it suffices to show
that  $\fSgx  {(H_\xi\*X)}\nprc_a \fSgx  {(G\*W)}$.

Since 
 $\rst{p_{_\xi}}{{}^*E_\xi}$ is realized by $\bar a_\xi$ in $\fSgx 
{(G\*U)}$, by construction, we see from (3) and the definition of $W$ that it
must also be realized by $\bar a_\xi$ in $\fSgx  {(G\*W)}$.  Also,   
\num 4{$\fSgx  {(H_\xi\*W)}\prc \fSgx 
{(H_\xi\*U)}$\quad,}
\endnum 
by Lemma 1(v).  Since $\rst{p_{_\xi}}{{}^*E_\xi}$ is not realized in $\fSgx 
{(H_\xi\*U)}$, by construction, we see from (4) that it is not realized in
$\fSgx  {(H_\xi\*W)}$.

We now transfer this situation to a strong type over 
$\fSgx  {(H_\xi\*X)}$.  Because
$\fW$ is countable and dense, and $X$ is infinite, there
is an  automorphism of $\fW$ taking $U_0$ into $X$ 
(this is the point of introducing $W$ to replace $U$). 
Extend it to an automorphism $\vartheta$ of $\fN' = \fSg {M'} {(W)}$.  Then
appropriate restrictions of $\vartheta$ are automorphisms of $\fSgx  {(G\*W)} $
and $\fSgx  {(H_\xi\*W)}$.  Therefore,  
\num 5 {$\rst{\vartheta(p_{_\xi})}{{}^*\vartheta(E_\xi)}$ is realized in $\fSgx 
{(G\*W)}$ (by $\vartheta(\bar a_\xi)$), but not in $\fSgx 
{(H_\xi\*W)}$\quad.}  \endnum

Suppose, now, for contradiction, that $\fSgx  {(H_\xi\*X)}\prc_a
\fSgx  {(G\*W)}$.  Then 
\num 6 {$\rst{\vartheta(p_{_\xi})}{{}^*\vartheta(E_\xi)}$ is
realized in $\fSgx  {(H_\xi\*X)}$\quad.}
\endnum  
On the other hand, this supposition also gives
\num 7 {$\fSgx  {(H_\xi\*X)}\prc
\fSgx  {(G\*W)}$\quad.}
\endnum 
Since $\fSgx  {(H_\xi\*U)}\prc
\fSgx  {(G\*U)}$, by construction, we get from (3) and (4) that  
\num 8 {$\fSgx  {(H_\xi\*W)}\prc \fSgx  {(G\*W)}$\quad.}
\endnum  
Combining (7) and (8), we arrive at $\fSgx  {(H_\xi\*X)}\prc\fSgx 
{(H_\xi\*W)}$.  In view of (6), this forces
$\rst{\vartheta(p_{_\xi})}{{}^*\vartheta(E_\xi)}$ to be realized
in $\fSgx  {(H_\xi\*W)}$, which contradicts (5).  This proves (2)

We now construct a $\xi<\omega_1$ and a countable $X\subseteq U$ 
that contradict (2). Since $\fB$ is strictly $\ka$-generated, by
Lemma 2, it is isomorphic to $\fF_\ka$.  Thus, there is a set $V$
of generators of $\fB$, of size $\ka$, that is totally indiscernible
in $\fB$.  Exactly as in the proof of Lemma 3, we construct 
increasing sequences, $\langle\,X_n:n\in\omega\,\rangle$ and 
$\langle\,Y_n:n\in\omega\,\rangle$, of countable subsets of $U$ and $V$ 
respectively, and an increasing sequence $\langle\,\rho_n:n<\omega\,\rangle$
of countable ordinals, such that 
$$\fSgx  {(H_{\rho_n}\*X_n)}\subseteq\fSgx  {(Y_n)}\subseteq 
\fSgx  {(H_{\rho_{n+1}}\*X_{n+1})}$$
for every $n<\omega$.  Set 
$$\xi=\sup\{\rho_n:n<\omega\}\quad,\quad 
X=\bigcup_{n\in\omega}X_n\quad,\quad\text{and}\quad
Y=\bigcup_{n<\omega}Y_n\quad.$$ 
Then 
$$\fSgx 
{(H_\xi\*X)}={\tsize\bigcup}_{n<\omega}\fSgx  {(H_{\rho_n}\*X_n)}=
{\tsize\bigcup}_{n<\omega}\fSgx  {(Y_n)}=\fSgx  {(Y)}\quad.$$

Since $\fSgx {(Y)}\prc_a\fB$ by Lemma 8 (with $K=\varnothing$) and the total
indiscernibility of $Y$ and $V$, we conclude that $\fSg M
{(H_\xi\*X)}\prc_a \fB$.  This is just the desired contradiction to (2).
\qed\enddemo

\proclaim{Lemma 10} 
For
every countable $H\subseteq L'-L$\comma there is a countable $K$\comma with
$H\subseteq K\subseteq L'\diff L$\comma such that\comma whenever $\fV$ is a dense
ordering and $\fN'=\EM V$\comma we have $$\fSg N {(K\*V)} \prc_a \fN\quad.$$
\endproclaim 
\demo{Proof}
The proof is almost identical to the proof of the second part of Lemma 4,
with ``$\prc$" replaced everywhere by ``$\prc_a$". One uses Lemma 9 in
place of the first part of the proof of Lemma 4, and Lemma 8 in place
of Lemma 1(v) and Facts 1 and 3.  Notice, for example, that, taking 
$K=L'\diff L$, we get essentially Fact 3 with ``$\prc$" replaced by
``$\prc_a$".  We leave the details to the reader.
\qed\enddemo

\proclaim{Lemma 11} $\fF_\lambda$ is
$a$-saturated for every $\lambda\geq\omega$.
\endproclaim
\demo{Proof}
By Shelah \cite{1987}, the proof of Theorem 2.1 on pp\. 285--287, $T$ has an
Ehrenfeucht-Mostowski model $\fN' = EM(\fW,\Delta)$ of power $>
\max\{(2^{|T|})^+,\ka\}$ such that $\fN = \rst {\fN'} L$ is $a$-saturated.  By
Lemma 10 there is a countable $K\subseteq L'\diff L$ such that $\fSg N
{(K\*W)}\prc_a \fN$.  Let $\fU$ be a dense submodel of $\fW$ of power $\ka$. 
Then  $\fSg N {(K\*U)}\prc_a\fSg N {(K\*W)}$, by Lemma 8.  Thus, we immediately
see that $\fSg N {(K\*U)}$, too, is $a$-saturated.  Since it is strictly
$\ka$-generated, by Lemma 2, it is isomorphic to $\fF_\ka$, by
$\ka$-uniqueness. Thus, $\fF_\ka$ is $a$-saturated. 

Next, we turn to $\fF_\omega$. Recall that
$\fF_\ka$ has a $\Phi$-basis $Z_\ka$ that extends a $\Phi$-basis $Z_\omega$ of
$\fF_\omega$\p. Hence, $\fF_\omega\prc\fF_\ka$ and even
$\fF_\omega\prc_a\fF_\ka$, by Lemmas 1(v)  and 8 (with $K=\varnothing$) and
the total indiscernibility of $Z_\ka$ and $Z_\omega$. Since $\fF_\ka$ is
$a$-saturated, by the previous lemma, so is $\fF_\omega$.  

Now consider any $\lambda>\omega$.  To show that 
$\fF_\lambda$ is $a$-saturated, consider an arbitrary 
strong typed $p$ based on a finite subset $E$ 
of $F_\lambda$\p.  Let $V$ a denumerably infinite 
subset of the $\Phi$-basis  $Z_\lambda$ of 
$\fF_\lambda$ that contains a set of generators for $E$.
Since $\fSg{} {(V)}$ is isomorphic to $\fF_\omega$\p, 
and the latter is $a$-saturated, $\pe$ is realized in
$\fSg{} {(V)}$. But  $\fSg{} {(V)}\prc\fF_\lambda$\p,
so $\pe$ is also realized in $\fF_\lambda$\p.
\qed\enddemo

\proclaim{Lemma 12} There is a cardinal $\lambda$ such that 
every model of $T$ of power $>\lambda$ is $a$-saturated.
\endproclaim
\demo{Proof} Let $\lambda$ be the Hanf number for omitting types in languages of
cardinality at most $|T|$.
 Suppose, for contradiction, that there  is a model $\fA$ of $T$ of power
 $>\max\{\lambda,\ka\}$ that is not $a$-saturated.  Then there is a strong type
$p = \pe$
 (over $\fA$) based on a finite subset $E$ of $A$ such that $p$ is omitted
 by $\fA$. 
Since there are at most $|T|$
many inequivalent formulas in $p$, by Shelah \cite{1990}, Chapter III, Lemma 2.2(2),
we may assume that $p$ is a (possibly incomplete) {\it type\/}---as
opposed to a strong type---over some subset $B$ of $A$ of cardinality at most
$|T|$.

 Let $L'$ be a language with built-in Skolem functions that extends $L$ and 
 that includes constants for the elements of $B$. Let $T'$ be a Skolem theory
 in $L'$ that extends the theory of $\lr\fA b_{b\in B}$. By  Morley's
 Omitting Types Theorem and the choice
 of $\lambda$, there is an
Ehrenfeucht-Mostowski model $\fN'$ of $T'$ over a
 dense ordering $\fU$ of indiscernibles of cardinality $\geq\lambda$ such 
 that $\fN'$ omits the type $p$ (see Morley \cite{1965}, Theorem 3.1, or
Chang-Keisler \cite{1973}, Exercise 7.2.4). Let $\fV$ be a dense submodel of
$\fU$ of 
 power $\ka$,
and set $\fM'=\fSg {N'} {\!(V)}$. Then $\fM'$ is the corresponding 
Ehrenfeucht-Mostowski model over $\fV$, and $\fM'\prc\fN'$, by Fact 3.
Obviously, $\fM'$, and hence also $\fM$,  omits $p$.

Let $H\subseteq L'-L$ be a finite set containing all the individual constant symbols for
elements of $E$.  By Lemma 4, there is a countable $K$ with 
$H\subseteq K\subseteq L'-L$ such that $\fSgx {(K\*V)}\prc\fM$. 
Thus, $E$ is a subset of the universe of $\fSgx {(K\*V)}$.  Since $p$, as a
strong type, is  based on a subset of the model $\fSgx {(K\*V)}$, every formula
$\varphi$ of $p$ is equivalent to a formula $\varphi'$ with parameters from 
$\Sgx {(K\*V)}$, by Shelah {\it ibid.\/}, Lemma
2.15(1).  Let $p'$ be the set of these formulas $\varphi'$.  Then $p'$
is a strong type based on $E$ that is equivalent to $p$.  Because $\fM$ omits
$p$, it also omits $p'$. Therefore, the elementary substructure $\fSgx
{(K\*V)}$ must
 omit $p'$.
This shows that $\fSgx {(K\*V)}$ is not
$a$-saturated.
  But that is impossible:
since $\fSgx {(K\*V)}$ is  strictly $\ka$-generated (by Lemma 2), it is
isomorphic to $\fF_\ka$\p,  and we have seen that $\fF_\ka$ is $a$-saturated.
\qed\enddemo

\proclaim{Lemma 13} $\fF_\lambda$ is $a$-prime over any $\Phi$-basis\comma
for $\lambda\geq\omega$. Consequently\comma $\fF_\omega$ is also 
$a$-prime over $\varnothing$.
\endproclaim
\demo{Proof}
Let $U$ be a $\Phi$-basis of $\fF_\lambda$\p.  Then $U$ generates any given
sequence  $\bar a$ from $\fF_\lambda$\p, so  $tp(\bar a,U)$ is
clearly atomic.  In particular,  $\fF_\lambda$ is $a$-atomic over $U$. 
Moreover, $\fF_\lambda$  is also $a$-saturated,  by Lemma 11, and there
are obviously no Morley sequences in $\fF_\lambda$ over $U$. By Shelah's second
characterization theorem for $a$-prime models, $\fF_\lambda$ is
$a$-prime over $U$  (see Shelah \cite {1990}, Chapter IV,
Definition 4.3 and Theorem 4.18).

Now an $a$-prime model over a countable set is also $a$-prime over
$\varnothing$ (see {\it ibid\.\/}).  Hence, $\fF_\omega$ is also $a$-prime over
$\varnothing$. \qed\enddemo

We now work inside of some monster model.  Let $p_0,\ldots,p_{n-1}$ be pairwise
orthogonal regular types, all based, say, on a finite set $E$. By realizing the 
strong type of $E$ over $\varnothing$ in $\fpr \varnothing$, and then passing to automorphic 
copies of $p_0,\ldots,p_{n-1}$\p, we may assume that $E$ is a subset of the universe of
$\fpr \varnothing$. For each $i<n$, let $I_i$ be an infinite Morley sequence built from 
$\rst{p_i}{^*E}$.
\proclaim{Lemma 14} $\fpr {E\cup\bigcup_{i<n}I_i}$ is strictly
$\lambda$-generated\comma  where $\lambda=|\bigcup_{i<n}I_i|$.
\endproclaim
\demo{Proof} The proof is by induction on $\lambda$.  If $\lambda=\omega$,
then, as mentioned in the preceding lemma,  $\fpr {E\cup\bigcup_{i<n}I_i}$ is $a$-prime over
$\varnothing$. But, we just saw that $\fF_\omega$ is
$a$-prime over $\varnothing$.  Hence, these two models are isomorphic, by the
uniqueness of $a$-prime models.  Since $\fF_\omega$ is  strictly
$\omega$-generated, by Lemma 2, so is $\fpr {E\cup\bigcup_{i<n}I_i}$\p. 

Now suppose that the lemma is true for all infinite $\mu<\lambda$.  Represent 
$\bigcup_{i<n}I_i$ as a strictly increasing sequence $\langle\,\bar a_\alpha:
\alpha<\lambda\,\rangle$ with the property that, for each $i<n$, the sequence
contains infinitely many members of $I_i$.  For each infinite $\beta<\lambda$, set
$\fB_\beta =\fpr {E\cup\{\bar a_\alpha:\alpha<\beta\}}$. We may arrange this so
that $\fB_\beta\prc\fB_\eta$ whenever $\beta\leq\eta$, and, at limit stages
$\delta$, that $\fB_\delta=\bigcup_{\beta<\delta}\fB_\beta$\p.  By the induction 
hypothesis, each $\fB_\beta$ is $|\beta|$-generated, since
$\omega \leq |\beta|<\lambda$. Therefore, $\fpr
{E\cup\bigcup_{i<n}I_i}=\bigcup_{\beta<\lambda}\fB_\beta$ is  generated by a set of
cardinality at most $\lambda\cdot\lambda=\lambda$.

To see that no set of smaller cardinality can generate it, we proceed by 
contradiction. Let $X$ be an infinite set of power $\nu<\lambda$ that
generates $\fpr {E\cup\bigcup_{i<n}I_i}$\p. Since $\fpr {E\cup\bigcup_{i<n}I_i}$ is 
$a$-atomic over $E\cup\bigcup_{i<n}I_i$\p, for each finite sequence $\bar b$ from $X$
there is a finite subset $I_{\bar b}\subseteq\bigcup_{i<n}I_i$ such that
$$stp(\bar b,E\cup I_{\bar b}) \vdash stp(\bar b,E\cup{\tsize
\bigcup}_{i<n}I_i)\quad.$$ Set $J=\bigcup\{I_{\bar b}: \bar b \text{ is a finite
sequence from } X\}$. Then $J$ has power at most $\nu$ and 
\num 1 {$stp(\bar b,E\cup J) \vdash stp(\bar b,E\cup{\tsize\bigcup}_{i<n}I_i)$
for every (finite) sequence $\bar b$ from $X$\quad.}
\endnum
 Since $\nu<\lambda$, there
is a sequence $\bar a=\langle\,a_0,\ldots,a_{k-1}\,\rangle$ that is in
$\bigcup_{i<n}I_i$, but not in $J$.  Thus, $\bar a$ is independent over $E\cup
J$.  This is readily seen to contradict (1).  For example, because $X$
generates $\fpr {E\cup\bigcup_{i<n}I_i}$\p, there is a sequence  $\bar b$
from $X$ and terms $\sigma_i(\bar x)$ of $L$, for $i<k$, such that
$\sigma_i[\bar b]=a_i$\p.  Now $$stp(\bar b,E\cup J) \vdash \bigwedge_{i<k}
\sigma_i(\bar x)=a_i\quad,$$ by (1), so there is a formula $\varphi(\bar x)$
in $L^*(E\cup J)$  (the set of formulas almost over $E\cup J$) such that 
$$\vdash \varphi(\bar x) \to \bigwedge_{i<k} \sigma_i(\bar x)=a_i\quad.$$
Thus, $\bar a$ is definable, over $L^*(E\cup J)$, since, e.g.,  $$\vdash
y=a_i \leftrightarrow \exists{\bar x}(\varphi(\bar x)\wedge  \sigma_i(\bar
x)=y)\quad \text{ for each } i<k\quad.$$ But then $\bar a$ is in the
algebraic closure of $E\cup J$, so it cannot be  independent over this set.
We have reached our contradiction. \qed\enddemo

\proclaim{Theorem 15} $T$ is categorical in every cardinality $>|T|$.
\endproclaim
\demo{Proof} 
We begin by proving that $T$ must be unidimensional.  Indeed, suppose for 
 contradiction that there are two orthogonal regular types, $p_1$ 
and $p_2$\p, over a finite subset $E$ of $\fpr \varnothing$. Working inside a 
monster model, for $i=1,2$, let $I_i$ be a Morley sequence of cardinality $\ka$
built from $\rst{p_i}{^*E}$. Then $\fpr{E\cup I_1}$ and $\fpr{E\cup I_1\cup I_2}$
are both strictly $\ka$-generated models, by the previous lemma, and 
hence isomorphic, by $\ka$-uniqueness. 
But this is impossible: the first model has a single 
dimension 
of cardinality $\ka$, the $p_1$-dimension, and all the other dimensions are
$\omega$;  the second model has exactly two dimensions of cardinality $\ka$, the
$p_1$- and the $p_2$-dimensions.  Thus, $T$ can only have one regular type, up to
equivalence. Let's denote this type by $p$.
\par
Now let $\lambda$ ($\geq|T|$)  be so large that all the models of $T$ of
power $> \lambda$ are $a$-saturated (see Lemma 12).  Then $T$ is 
$\lambda^+$-categorical. Indeed, if $\fM$ and $\fN$ are models of power 
$\lambda^+$, then they both are $a$-saturated, by choice of $\lambda$, 
and have $p$-dimension $\lambda^+$, by unidimensionality.  Hence, 
they must be isomorphic, by Shelah \cite{1990}, Chapter V, Theorem 2.10.  This
shows that $T$ is categorical in some power $>|T|$.  By the results of Shelah
\cite{1974}, it is categorical in every power $>|T|$.
 \qed\enddemo 

\proclaim{Theorem 16} $T$ is $\lambda$-unique for every uncountable $\lambda$.
\endproclaim
\demo{Proof} If $T$ is definitionally equivalent to a countable theory, and is
categorical in every infinite power, i.e., if $T$ is a totally categorical
theory, then the theorem is obvious.  Suppose, now, that it is not totally
categorical. We shall use Theorems 1--3, including the proof of Theorem 2, from
Laskowski \cite{1988}. According to these, $T_\infty$ must have a minimal
prime model $\fM_0$.  Moreover, there is a type $p$ based on $M_0$ with the
following properties: 

\medskip \smallskip\flushpar
 \halign{#\hfill&\quad#\hfill\cr
(1)&If $\fN$ is any elementary extension of $\fM_0$, then any two maximal
Morley sequences in\cr
\ &  $\fN$ built from $\rst p {M_0}$ must have the same
cardinality.  This is called the $p$-dimension\cr
\ &  of $\fN$.\cr}  
\medskip\flushpar
 \halign{#\hfill&\quad#\hfill\cr
(2)&If $\fN$ and $\fN'$ are elementary extensions of $\fM_0$, with the same
$p$-dimension, then  \cr 
\ & there is an isomorphism of $\fN$ onto $\fN'$ that fixes ${M_0}$ .\cr}  
\medskip\flushpar
 \halign{#\hfill&\quad#\hfill\cr
(3)&If $\fM_0 \prc \fN \prec\fN'$, then $\rst p N$ is realized in $\fN'$.
\cr}  
\medskip\smallskip\flushpar
To simplify notation, we shall assume that $p$ is a $1$-type.

 Since $T_\infty$ is a universal-existential theory categorical
in power, it is model complete, by Lindstr\"om's theorem (see, e.g.,
Chang-Keisler \cite{1973}, Theorem 3.1.12).  

Fix $\lambda > \omega$, and let $\fM$ be a strictly $\lambda$-generated model of
$T$, say $X$ is a generating set of power $\lambda$.  Our goal is to prove that
$\fM$ has $p$-dimension $\lambda$.  Let $W$ be a countably infinite subset of
$X$.  Then $\fSgx {(W)} \prc \fM$, by model completeness.  Without loss of
generality, we may suppose that $\fM_0$ is an elementary submodel of $\fSgx
{(W)}$.  We now define a strictly increasing, continuous sequence,
$\langle Y_\xi:\xi<\lambda\rangle$, of subsets of $X$ of size $<\lambda$, such
that   \num 4 {$\fM_0 \prc \fSgx {(Y_\xi)} \prec \fSgx
{(Y_{\xi+1})}\prc \fM$\quad for $\xi<\lambda$\quad.}
\endnum
Indeed, take $Y_0 = W$.  If $Y_\xi$ has been defined, then $\fSgx {(Y_\xi)}\neq
\fM$, since $\fM$ is strictly $\lambda$-generated and $Y_\xi$ has power
$<\lambda$.  Hence, there is an element $u_\xi$ in $X$ that is not generated by
$Y_\xi$.  We set $Y_{\xi + 1}=Y_\xi\cup\{u_\xi\}$.  Property (4) follows by our
choice of $Y_0$ and $u_\xi$, for each $\xi$, and by model completeness.

For each $\xi <\lambda$, the type $\rst p {\Sgx {(Y_\xi)}}$ is realized in
$\fSgx {(Y_{\xi + 1})}$, by (3) and (4).  Therefore, using our strictly increasing
 chain $\langle \fSgx {(Y_\xi)}: \xi<\lambda\rangle$ of elementary
substructures of $\fM$, we can build from $\rst p {M_0}$ a Morley sequence in
$\fM$ of length at least $\lambda$.  Thus, $\fM$ has $p$-dimension at least
$\lambda$.

Suppose now that $I$ is a maximal Morley sequence in $\fM$ built from $\rst p
{M_0}$.  Then $\fSgx {(M_0\cup I)}= \fM$.  For otherwise we would have 
$\fSgx {(M_0\cup I)}\prec \fM$, by model completeness.  Hence, we could realize
$\rst p {\Sgx {(M_0\cup I)}}$ in $\fM$, by (3), and thus extend $I$ to a
larger Morley sequence in $\fM$ built from $\rst p {M_0}$, contradicting its
maximality.  Since $X$ also generates $\fM$, and has cardinality
$\lambda>\omega$, there is a subset $J$ of $I$ of cardinality $\leq \lambda$
such that $M_0\cup J$ generates $X$, and hence also $\fM$.  In particular,
$M_0\cup J$ generates $I$.  But then $J=I$, since any element in $I\diff J$
would be independent over $M_0\cup J$, and hence could not be generated by
this set.  We conclude that $|I|\leq\lambda$.  In other words, $\fM$
has $p$-dimension at most $\lambda$, and hence exactly $\lambda$.

We have shown: 
\num 6 {Any strictly $\lambda$-generated model of $T$ extending $\fM_0$ has
$p$-dimension $\lambda$\quad.}
\endnum
Now let $\fM$ and $\fN$ be any two strictly $\lambda$-generated models of $T$. 
By passing to isomorphic copies, we may assume that $\fM_0$ is an elementary
substructure of each.  Hence, by (2) and (6), $\fM$ and $\fN$ are isomorphic over
$M_0$.  This shows that $T$ is $\lambda$-unique.\qed
\enddemo 

By the preceding theorem, the noncountably generated models of $T$ are, up to
isomorphisms, precisely the structures $\fF_\lambda$, for $\lambda >\omega$.

\enddocument